\documentclass[12pt,oneside]{article}
\newcommand{\Path}{}
\usepackage{\Path Paper_style_2}
\usepackage{\Path Keywords_WB}
\usepackage{\Path Environments}
\usepackage{\Path Category}
\newcommand{\Diagrams}{}
\usepackage{tikz-cd}
\usepackage{\Path tikzcd2}

\begin{document}
	
\title{Cells, convexity and contractibility in general categories}
\author{Suddhasattwa Das\footnotemark[1]}
\footnotetext[1]{Department of Mathematics and Statistics, Texas Tech University, Texas, USA}
\date{\today}
\maketitle

\begin{abstract} 
	The two pillars of Algebraic topology - homology and homotopy theory rely on the availability of basic building blocks called cells. Cells take the form of simplexes, and have properties such as faces, sub-cells, convexity and contractibility. The first two cells, namely the line and point lead to the concept of homotopy. The collection of maps from the cells and the redundancies among them determine the homology of objects. This article presents a procedure by which such cells can be built in general categories satisfying some simple axioms. The cells satisfy the categorical analogs of convexity and contractibility. This enables a cellular theory for the general category, carrying notions of homotopy, homology, cellular approximation and homotopy equivalence which are mutually compatible in the same way as in the familiar context of Topology.
\end{abstract}

\section{Introduction} \label{sec:intro}

Many disciplines in Mathematics have an underlying principle of building objects using a limited collection of simpler objects, so that the composite object is a sum of its parts and displays a richer set of properties. The best example of this approach is the theory of cell and $CW$-complexes in Topology. There has been a number of independent efforts in other disciplines made in the same spirit. The concept of cellular algebras \cite{GrahamLehrer1996cell} in Representation theory was proposed to study Kazhdan-Lusztig representations and Hecke algebras. This concept was later lifted to abstract category theory to create the notion of \textit{cellular categories} \cite{Westbury2009cell, AndersenEtAl2015cell, Tubbenhauer2023web, LiebermanEtAl2019cell}. Another important example is the use of \textit{cellular resolutions} and \emph{algebraic mapping cone}s for a cellular theory of chain complexes \cite{MedinaMardones2023cell, FilakovskyVokrinek2023cell, EngstromNoren2012cell, OlteanuWelker2016cell}. 

However, these ideas have been mutually distant and independent. There is a need to find conditions on a general category $\calC$ that would allow a \textit{cellular theory}. This means a sequence of objects called cells, which play the role of building blocks. Such cells should be equipped with well defined components called faces, edges and corners. A natural candidate for gluing would be the categorical operation of pushout. The most desirable property is that these cells display appropriately abstracted notions of convexity and contractibility. The cells themselves should be contractible so that all holes present in a cellular complex is a consequence of the arrangement and not from an individual cell. The cells should also be convex in a suitable sense, so that one can perform an analog of barycentric subdivision for closer simplicial approximations.

\begin{figure} [!t]
\centering
\begin{tikzpicture}[scale=1.0, transform shape, framed, background rectangle/.style={double, ultra thick, draw=gray, rounded corners}]
	\input{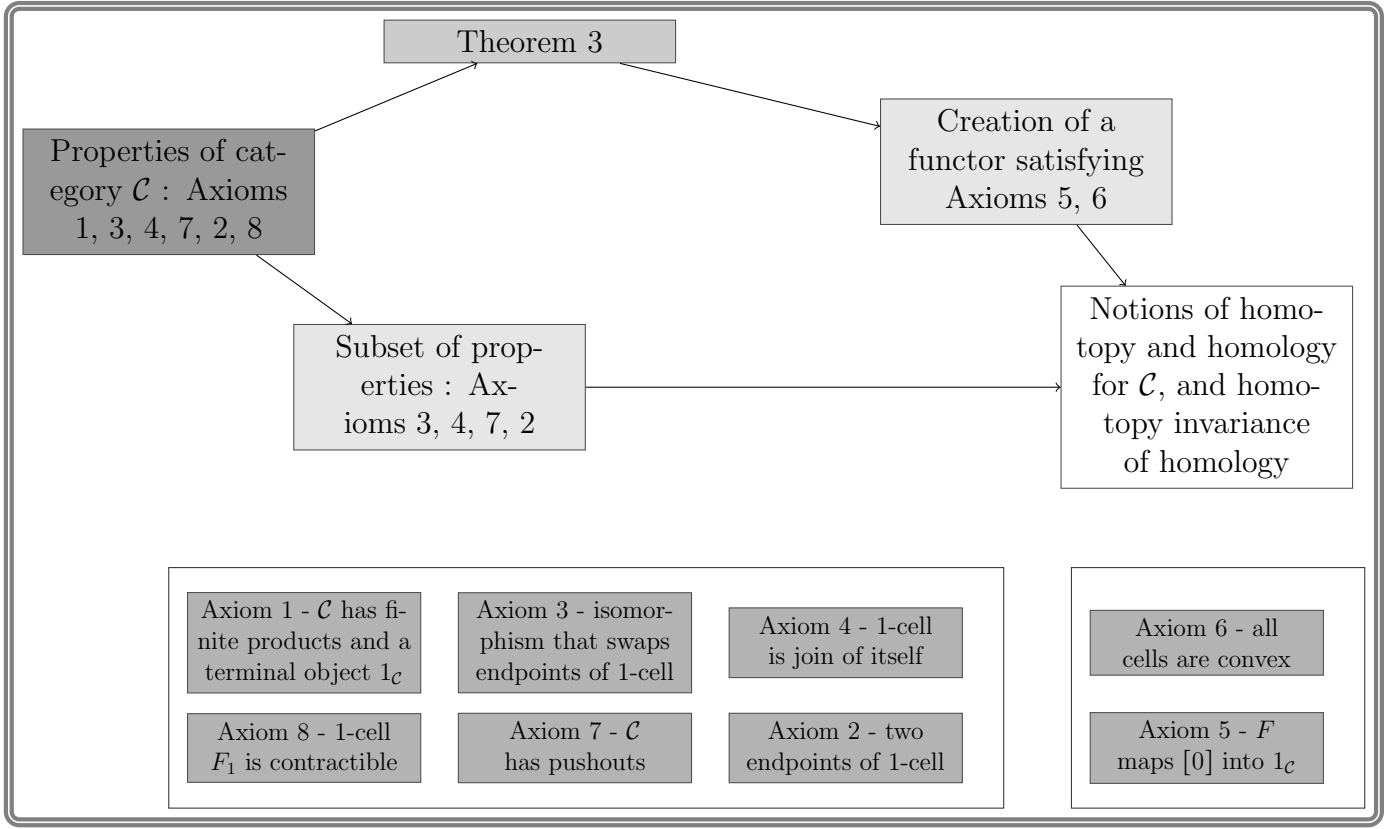}  
\end{tikzpicture}
\caption{Outline of the theory. The paper presents a simple axiomatic approach to both homotopy and homology which guarantee that the homology preserves homotopy invariance. }
\label{fig:outline:4}
\end{figure}

\begin{figure}[!t]
\centering
\begin{tikzpicture}[scale=1.0, transform shape, framed, background rectangle/.style={double, ultra thick, draw=gray, rounded corners}]
	\input{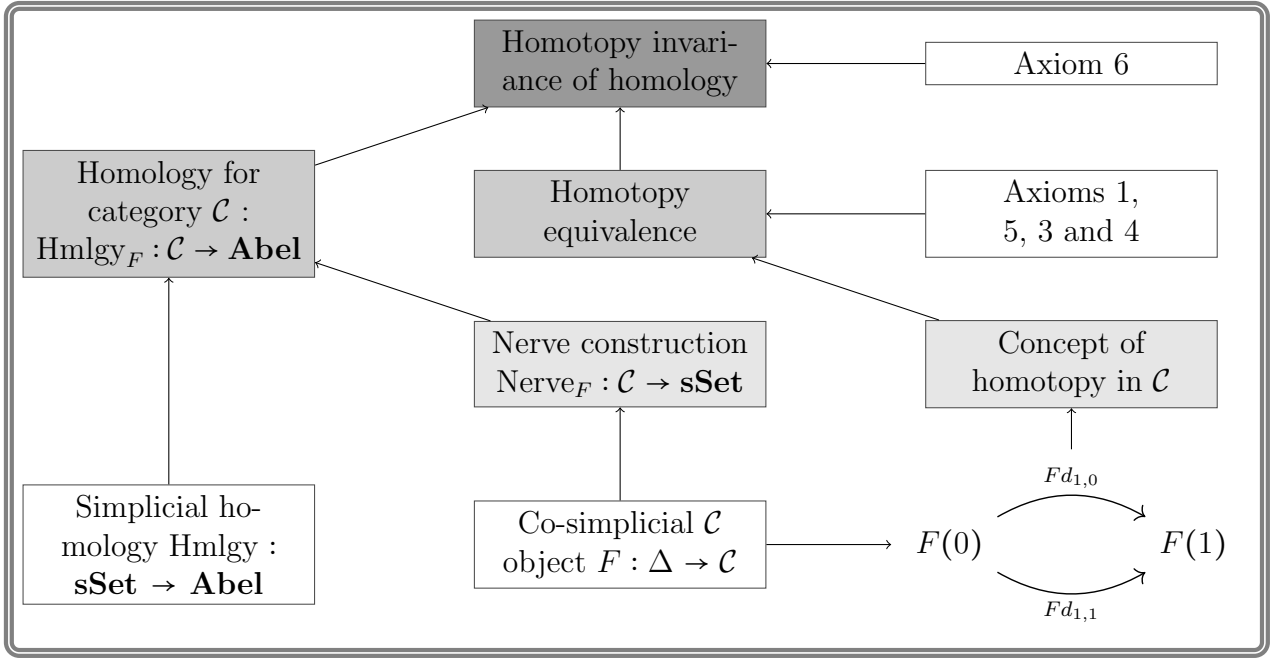}
\end{tikzpicture}
\caption{Logical dependencies of homotopy, homology and various categorical axioms. The chart above presents various assumptions and concepts, along with their logical dependence. The independent notions are in white boxes. This includes simplicial homology, a means of converting the combinatorial structure of a simplicial set into an Abelian group. Homology is the confluence of simplicial homology and the nerve construction. The nerve construction \eqref{eqn:def:Nerve} is an outcome of an arbitrary functor $F : \Delta \to \calC$ as in \eqref{eqn:functorF}. If one isolates the action of $F$ on just the first two objects of the simplex category $\Delta$, then one obtains a generalized notion of homotopy. This attains all the usual properties of homotopy equivalence under further axioms. See Figure \ref{fig:hmtpy_outline} for a detailed outline of this property. }
\label{fig:outline:2}
\end{figure}

Recent findings in a category theory \cite{Das2024hmlgy} have revealed that the topological notions of homotopy, homology, contractibility, convexity, and homotopy invariance of homology, can be created in a general category $\calC$ if one has a co-simplicial object, which is a functor $F : \Delta \to \calC$ from the combinatorial simplex category $\Delta$. The codomain $\calC$ will be called the \emph{context}. The required axioms on $\calC$ are displayed in Figures \ref{fig:outline:4} and \ref{fig:outline:2}.

The axiomatization in \cite{Das2024hmlgy} provides a simple common ground for both homology and homotopy, without using the heavy machinery of \textit{higher categories}.  The key idea was to distill the categorical or structural nature from some fundamental notions of topology, such as homotopy, contractibility and convexity. Instead of relying on their definitions based on points, sets, maps and convex sums,  these can be re-interpreted by their relational properties with other objects in the collection $\Topo$ of topological spaces and continuous maps. This re-interpretation allows a generalization to general category $\calC$, by laying a foundation through set of axioms on the category $\calC$ and co-simplicial object $F$. 

Some of the axioms on $F$, particularly the axiom of convexity, is hard to verify in arbitrary categories. Convexity is an intuitive concept in topological or linear spaces, but more technical in abstract categories. Our usual notion of convexity involves a topological space as well as its embedding in a linear space. Two topologically isomorphic objects may be simultaneously convex and non-convex. Convexity is presented in \cite{Das2024hmlgy} as a property extrinsic to an object, and independent of notions of embedding. 

The goal of this article is to present a simple set of conditions under which these axioms are realized. These axioms, summarized in Figure \ref{fig:outline:4}, are entirely on $\calC$ and do not involve a functor $F$. We start with two basic building blocks denoted $F_0$ and $F_1$. These are objects of $\calC$ and $F_0$ coincides with the terminal element $1$ of $\calC$. The object $F_1$ plays the role of abstract interval or line segment.

The line of thought presented in this article may be viewed from two ends. On one side, it presents a \emph{synthetic} approach to Topology, in which various topological properties and theorems are rediscovered by their relational properties, rather than by principles of point-set topology. Such a synthetic approach towards Topology has been found to be rewarding and elegant \cite[e.g.]{ClementinoTholen1997sep, GiuliLapal2009, GiuliTholen2000open}. On the other side, this work can be viewed as a \emph{topologification} of a general category. Topologification is the transfer of the machinery of Topology to other categories, by formulating the appropriate axiomatic machinery. While \textit{sheaf theory} \cite{rosiak2022sheaf} and \emph{topos theory} \cite{Flori2013Topos, Leinster2010topos} remain some prime examples, it has also been done in other arenas such as dynamical systems \cite{DasSuda2024recon} and Analysis \cite{LowenLowen1987conv, Yildiz2015uniform, KhadimQasim2022uniform, ClementinoEtAl2025lax, EscardoStreicher2016seq, GutierresHofmann2008seq}.

A notable body of work carrying the same objectives, was the work of Voedovsky \textit{et. al.}  \cite{MorelVoevodsky1991A1, voevodsky1992hmlgy, voevodsky1996hmlgy} for \textit{model categories}. Their idea was a formalization of the intuition that the category of algebraic varieties is similar to the category $\Topo$ of topological spaces. They go on to define a model category
structure on the category of simplicial sheaves on a site with interval. 
Model categories \cite{Quillen2006hmtpy, Hovey2013quil} are an abstraction of homotopy theory of topological spaces and chain complexes. A slight relaxation of the axiomatic framework of model categories was \emph{homotopical categories} \cite[e.g.]{DwyerEtAl2005, Riehl2014homotopy}, with main emphasis laid on the so called \textit{2-out-of-6 property} for weak equivalences—was formally established to handle diagram categories lacking compatible model structures. Voedevosky \textit{et. al.}'s construction \cite{MorelVoevodsky1991A1} works only in certain sheaf or model categories, whereas the present article provides a specific recipe to meet the general framework of \cite{Das2024hmlgy}. The main distinguishing feature is attaining a categorical analog of \textit{convexity} for the cells which turns out to be necessary for homology to be homotopy equivalence. We shall explore some connections of this notion of convexity with \emph{decalage} \cite{Stevenson2011decalage, Stevenson2012decalage, Stevenson2016decalage, Lurie2006topos}. 

One of the key ideas to be presented is that of a \emph{wedge-construction} $\Wedge$. It is a categorical abstraction of the process by which any shape can be converted into a cylinder and then into a cone. Once the functorial properties of the wedge construction are established, it is used to recursively create a sequence of objects $F_{n+1} = \Wedge(F_n)$ for $n=1,2,3,\ldots$. The object $F_n$ thus created will be the categorical equivalent of a geometric simplex. These objects will be called \emph{cells}. The wedge construction is a categorical generalization of the classical topological mapping cone \cite[e.g.]{massey1991basic, Spanier2012algebraic, Hatcher2002algtopo}. It is a simple version of other extensions such as the \textit{algebraic mapping cone} \cite[e.g.]{Weibel1994Hmlgy} on chain complexes, and as \emph{distinguished triangles} \cite[e.g.]{neeman2001triangulated, KashiwaraSchapira2006cat}. The simpler version is sufficient for our purpose, and its functorial properties are investigated in depth.

The cells so constructed will be of fundamental importance in the study of $\calC$. They enable the creation of a \emph{cellular theory} for $\calC$. The wedge construction naturally encodes a notion of \textit{faces} and \textit{boundaries} on these objects. For every pair of indices $m<n$, the $m$-th cell $F_m$ can be embedded in $^{n+1}C_{m+1}$ ways  into the $n$-th cell $F_n$, analogous to how the $n$-th convex simplex has exactly $^{n+1}C_{m+1}$-many $m$-sub simplices. These properties allow cells to be attached together to create compound objects, similar to cell complexes in topology. They will be shown to have the categorical equivalent of the topological properties of contractible and convex.

The main result (Theorem \ref{thm:recipe}~(i)) is to show that these sequence of objects $\SetDef{F_n}{n\in \num_0}$ along with their boundary and face morphisms, assemble into a co-simplicial object $F: \Delta \to \calC$. Any such functorial object generates a homology for the category via the associated nerve functor. This coincides with the familiar notion of singular homology when $\calC = \Topo$. In the general case, it serves as an algebraic characterization of objects of $\calC$, assigning a sequence of Abelian groups to each object, and group homomorphisms for every morphism in way that preserves compositionality.

The second main result (Theorem \ref{thm:recipe}~(iii)) to be presented is that this assembly also satisfies an Axiom of convexity \ref{A:convex}. This involves a set of identities between the faces of the cells, which allows an easy computation of homology. Convexity is the key property that makes homology an invariant of homotopy, fact presented in detail in \cite{Das2024hmlgy} and outlined in Figure \ref{fig:outline:1}.

Overall, the analysis presented provides a deconstruction of the familiar notions of homotopy and homology into components which are not specific to the context of topology but are relational. This allows for a broader interpretation of these notions and constructions. A particularly interesting prospect is the development of a cellular theory for general categories. Section \ref{sec:conclus} presents some idea on how these offers the possibility of cellular approximation and triangulation.

\section{Categorical origin of homology and homotopy} \label{sen:hmlgy_hmtpy}

We start with the most basic axiom on our context category $\calC$.

\begin{Axiom} \label{A:C}
The category $\calC$ has a terminal object $1_{\calC}$, and has finite products.
\end{Axiom}

A terminal object is a distinguished object $1_{\calC}$ such that for any object $x$ of $\calC$, there is a unique arrow from $x$ to $1_{\calC}$. In $\Topo$ and $\Delta$ the terminal objects are respectively $\star$ the 1-point space, and $0$. Categorical products is a generalization of the concept of Cartesian products. Given two objects $a,b$, their product if it exists, is an object denoted as $a \times b$ along with morphisms called \emph{projections} $a \xleftarrow{\proj_1} a\times b \xrightarrow{\proj_2} b$ such that for any third object $c$ and any pair of morphisms $c \xrightarrow{f} a$ and $c \xrightarrow{g} b$, there is a unique morphism $c \xrightarrow{ \langle f,g \rangle } a \times b$ through which both $f,g$ factors as shown below : 
\[\begin{tikzcd}
&  & c \arrow[lld, "f"'] \arrow[rrd, "g"] \arrow[d, "{\langle f, g \rangle}"] &  &   \\
a &  & a \times b \arrow[ll, "\proj_1"] \arrow[rr, "\proj_2"']                  &  & b
\end{tikzcd}\]
In $\Topo$, the categorical product of two spaces is the same as their Cartesian products.

\begin{Axiom} \label{A:brace}
There is an object $F_1$ in $\calC$ along with morphisms $\tilde{d}_{1,0}, \tilde{d}_{1,1} : 1_{\calC} \to F_1$ .
\end{Axiom} 

The object $F_1$ plays the role of an interval in $\Topo$. The two morphisms in Axiom \ref{A:brace} play the role of the two endpoints of the interval. In a categorical approach such as ours, one does not have sets and points. Instead the role of points is abstracted by morphisms from the terminal object $1$.

\begin{figure} [!t]
\centering
\begin{tikzpicture}[scale=1.0, transform shape, framed, background rectangle/.style={double, ultra thick, draw=gray, rounded corners}]
	\input{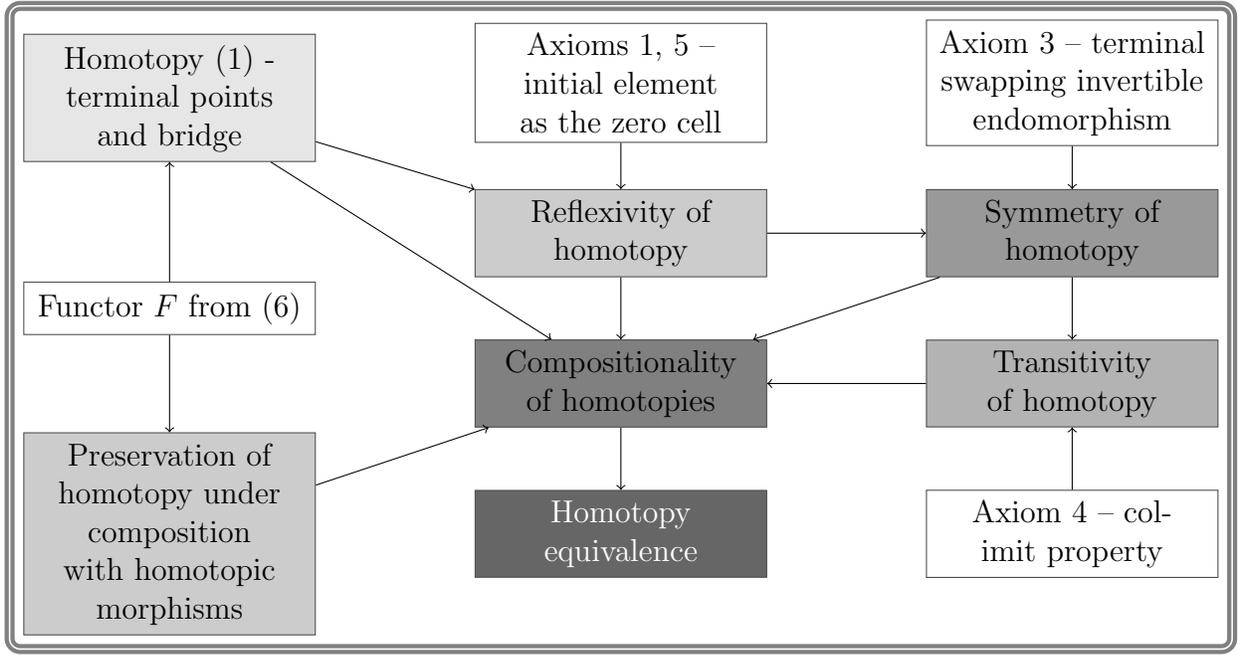}
\end{tikzpicture}
\caption{Construction of homotopy.}
\label{fig:hmtpy_outline}
\end{figure}

\paragraph{Homotopy} We briefly review how Axioms \ref{A:C} and \ref{A:brace} enable the concept of homotopy with all its recognizable properties within the general category $\calC$. We follow the outline in Figure \ref{fig:hmtpy_outline}, presented in \cite{Das2024hmlgy}. Given two objects $X,Y\in ob(\calC)$, we say that two morphisms $f,g:X\to Y$ are \emph{homotopic} if there is a morphism $H : X\times F_1 \to Y$ such that the following commutation holds :
\begin{equation} \label{eqn:def:hmtpy:2}
\begin{tikzcd} 
	X\times 1_{\calC} \arrow[Holud]{rrr}{ X \times F\paran{ d_{1,0} } } &&& X\times F(1) \arrow[Shobuj]{d}{H} &&& X\times 1_{\calC} \arrow[Akashi]{lll}[swap]{ X \times F\paran{ d_{1,1} } } \\ 
	X \arrow[Holud]{rrr}[swap]{f} \arrow[Holud]{u}{\cong} &&& Y &&& X \arrow[Akashi]{u}[swap]{\cong} \arrow[Akashi]{lll}{g}
\end{tikzcd}
\end{equation}
Thus the $1$-cell $F_1$ serves as a bridge between its two endpoints. A homotopy between $f$ and $g$ is a morphism $H$ from the cylinder $X\times F_1$ into $Y$, such that when restricted to its end faces, one gets $f$ and $g$ respectively. Henceforth, we shall use $F_0$ to denote the terminal object $1$. One of the most important notions associated to homotopy is being \textit{contractible}.

\paragraph{Contractibility} Similar to classic Homotopy theory, an object $X$ of $\calC$ will be called \emph{contractible} if its identity morphism $X$ is homotopic to some constant endomorphism. A \emph{constant endomorphism} is a composite of morphisms of the form $X \xrightarrow{!} 1_{\calC} \xrightarrow{x} X$. Here $x$ is an element of $X$ and the composite morphism can be interpreted to be constant of value $x$. Applying the commutative definition \eqref{eqn:def:hmtpy:2} of homotopy to contractibility, we get the diagram on the left : 
\begin{equation} \label{eqn:conractible:3}
\begin{tikzcd}
	X\times F_0 \arrow[ rr, "X \times \tilde{d}_{1,0}" ] && X \times F_1 \arrow[d, "\Contraction_X"] && X\times F_0 \arrow[ ll, "X \times \tilde{d}_{1,1}"' ] \\
	X \arrow[u, "="] \arrow[rr, "="] && X & F_0 \arrow[l, "x"] & X \arrow[u, "="'] \arrow[l, "!_X"]
\end{tikzcd} \quad \Leftrightarrow\quad 
\begin{tikzcd}
	X \times F_1 \arrow[dr, "\Contraction_X"] && X\times F_0 \arrow[ ll, "X \times \tilde{d}_{1,1}"' ] \arrow[d, "!"] \\
	X\times F_0 \arrow[ u, "X \times \tilde{d}_{1,0}" ]\arrow[r, "="'] &  X & F_0 \arrow[l, "x"] & 
\end{tikzcd}
\end{equation}
The diagram on the right is a simplified version of the diagram on the left. Equation \eqref{eqn:conractible:3} is an equivalent definition of contractibility. 

The following axiom is needed to make homotopy a symmetric relation 

\begin{Axiom} \label{A:swap}
There is an isomorphism $\text{swap} : F_1 \to F_1$ such that the following commutation holds with boundary maps :
\begin{equation} \label{eqn:def:hmtpy:1}
	\begin{tikzcd}
		& & F_1 \arrow[d, "\text{swap}", bend left=49] \\
		1_{\calC} \arrow[rr, "\tilde{d}_{1,1}"'] \arrow[rru, " \tilde{d}_{1,0}", bend left] & & F_1 \arrow[u, "\text{swap}^{-1}", bend left=49]
	\end{tikzcd}
\end{equation}
\end{Axiom}

The following axiom is required to make homotopy a transitive property.

\begin{Axiom} \label{A:F1_join}
The diagram 
\begin{equation} \label{eqn:F1_join:1}
	\begin{tikzcd} [column sep = large]
		F_1 & F_0 \arrow[l, "\tilde{d}_{1,0}"'] \arrow[r, " \tilde{d}_{1,1}"] & F_1
	\end{tikzcd}
\end{equation}
has a colimit, which is $F_1$ itself.
\end{Axiom}

The colimit diagram is a pushout square as shown below
\begin{equation} \label{eqn:F1_join:2}
\begin{tikzcd}
	F_0 \arrow[dotted, Shobuj, drrr, "\text{mid}"] \arrow[rrr, "\tilde{d_{1,0}}"] \arrow[d, "\tilde{d_{1,1}}"'] &&& F_1 \arrow[dotted, Shobuj, d, "\text{left}"] \\
	F_1 \arrow[dotted, Shobuj, rrr, "\text{right}"'] &&& F_1
\end{tikzcd}
\end{equation}
These morphisms capture the intuitive notion that a line segment joined end-to-end is a homeomorphic line segment. The two pieces being joined form the left and right halves respectively. The common point of gluing remains in the middle of the joined segment. The utility of the Axioms developed so far is summarized below :

\begin{proposition} [Generalized homotopy] \label{prop:hmtpy}
\cite{Das2024hmlgy} Under Axioms \ref{A:C}, \ref{A:brace} there is a notion of homotopy between morphisms of $\calC$, as defined in \eqref{eqn:def:hmtpy:2}. Further, this relation is symmetric under Axiom \ref{A:swap} and transitive under Axiom \ref{A:F1_join}.
\end{proposition}

Proposition \ref{prop:hmtpy} helps to capture many of intuitive features of homotopy in the abstract category $\calC$. Here is an example involving the important notion of contractibility.

\begin{lemma} [Contractible objects] \label{lem:cntrctbl:1}
Suppose Axioms \ref{A:C}, \ref{A:brace} and \ref{A:F1_join} hold. Then the product of two contractible objects is contractible.
\end{lemma} 

To see why, take two contractible objects $A,B$ which can be contracted to points $a,b$ respectively. Then observe the following commutative diagram.
\[\begin{tikzcd}
A \times B \arrow[drr, Akashi, shift right=1pt, "a \times B"'] \arrow[drr, Holud, shift left=1pt, "A \times B"] \arrow[rrrr, Akashi, dashed, shift right=1pt, "a \times b"'] \arrow[rrrr, dashed, Holud, shift left=1pt, "A \times B"] &&&& A \times B \\
&& A \times B \arrow[urr, Akashi, shift right=1pt, "A \times b"'] \arrow[urr, Holud, shift left=1pt, "A \times B"]
\end{tikzcd}\]
Each parallel pair of solid blue and yellow arrows are homotopic by assumption. Axiom \ref{A:F1_join} enables transitivity. This means that their compositions are homotopic. They are precisely the horizontal dashed pair of arrows.

We now turn our attention to a fixed combinatorial category essential of all our constructions.

\paragraph{Simplex} The simplex category $\Delta$ has as objects the non-negative integers $\num_0 :- \braces{0, 1,2,\ldots }$. Each such integer $n$ is meant to represent the ordered sets $[n] := \{0, \ldots, n\}$. The morphisms are given by
\[ \Hom_{\Delta} (m;n) := \braces{ \mbox{ Order preserving maps } \phi:[m] \to [n] } . \]
The simplex category $\Delta$ is the most concise way of encoding combinatorial structures in other categories \cite{Riehl2011sset, Steiner2007omega, Grodal2002higher, friedman2008sSet}. It finds applications is many different branches of mathematics \cite[e.g.]{Meush1981cat, BubenikScott2014prstnt} which involves enumeration. Recall that the category $\sSet$ of \textit{simplicial sets} is the functor category $\Functor{\Delta^{op}}{\SetCat}$. This is a category with purely combinatorial content, yet is essential to most categorical studies. Similarly, a functor from $\Functor{\Delta}{\SetCat}$ is called a \textit{co-simplicial set}. More generally, a functor of the form
\begin{equation} \label{eqn:functorF}
F : \Delta \to \calC , 
\end{equation}
is called a \textit{co-simplicial $\calC$-object}. It is a sequence of objects in $\calC$ along with a record of morphisms between them which observe the combinatorial rules within $\Delta$. Co-simplicial objects have been used in many different ways to obtain useful nerve functors (see \eqref{eqn:def:Nerve}) and homology theories, such as \emph{motivic theory} \cite{SuslinVoevodsky1997singular}, for Dold-Kan correspondences \cite{dold1958homology, Kan1958css} and polynomial differential forms \cite{Sullivan1977infntsml}. We shall assume

\begin{Axiom} \label{A:1_0_cell}
The functor $F$ from \eqref{eqn:functorF} maps $[0]$ into the terminal element $1$ from Axiom \ref{A:C}.
\end{Axiom}

Axiom \ref{A:1_0_cell} is intuitive and turns out to be essential to establish a connection between homology and homotopy. Homology originates from the \emph{nerve} of the functor $F$, as explained next.

\paragraph{Nerve construction}  The general functor $F$ induces a special functor called the \emph{nerve} of the category $\calC$ \cite{Riehl2011sset, GoerssJardine2009} as shown below : 
\begin{equation} \label{eqn:def:Nerve}
\begin{split}
	& \Nerve = \Nerve_F : \calC \to \sSet, \\ 
	& \forall c\in ob(\calC) \;:\; \Nerve(c) := \Hom \paran{ F\cdot; c } : \Delta^{op} \to \SetCat .
\end{split}
\end{equation}
The functor $\Nerve_F$ assigns a simplicial set to each object of $\calC$. Thus it assigns a purely combinatorial identity to each object in $\calC$. This functor is the key ingredient to creating a notion of Homology. For the example of \eqref{eqn:def:StndrdTopoSmplx}, $\Nerve_F(X)$ is the collection of all possible continuous mappings of the $n$-dimensional simplex into the topological object $X$. Thus $\Nerve_F(X)$ becomes a ledger for the topological ``content" of $X$, and the entries of the ledger are all possible embeddings of $n$-dimensional disks. 

\paragraph{Base maps} Let $\tilde{\Delta}$ be the subcategory of $\Delta$ generated only by the face maps $\SetDef{ d_{n,i} }{ 0 \leq i \leq n }$. Let $\tilde{\Delta} \xrightarrow{i} \Delta$ denote the obvious inclusion of categories. The simplified simplex $\tilde{\Delta}$ has an in-built $\text{Shift}$ functor : 
\begin{equation} \label{eqn:def:Shift}
\begin{tikzcd} \tilde{\Delta} \arrow[rr, "\text{Shift}"] && \tilde{\Delta}  \end{tikzcd}, \quad 
\begin{tikzcd} \Holud{n} \arrow[d, "d_{n,i}"'] \\ \akashi{n+1} \end{tikzcd} 
\begin{tikzcd} {} \arrow[r, mapsto] & {} \end{tikzcd}
\begin{tikzcd} \Holud{n+1} \arrow[d, "d_{n+1,i+1}"] \\ \akashi{n+2}	\end{tikzcd}
\end{equation}
A key observation is the existence of a natural transformation $\Base$ as shown below :
\begin{equation} \label{eqn:def:Base}
\begin{tikzcd}
	& \Delta^{op} \arrow{rrr}{\Nerve_{F}(X)} &&& \SetCat \\
	\tilde{\Delta}^{op} \arrow[dashed]{urrrr}[name=n1]{} \arrow[d, "\text{Shift}"'] \arrow[dashed]{drrrr}[name=n2]{} \arrow[d, "\text{Shift}"'] \arrow[ur, "\iota"] \\ 
	\tilde{\Delta}^{op} \arrow[r, "\iota"] & \Delta^{op} \arrow{rrr}[swap]{\Nerve_{F}(X)} &&& \SetCat
	\arrow[shorten <=1pt, shorten >=1pt, Rightarrow, to path={(n2) to[out=90,in=-90] node[xshift=20pt]{$\Base$} (n1)} ]{  }
\end{tikzcd} , \quad 
\begin{tikzcd}
	\Hom_{\calC} \paran{ F(n+1); X } \arrow[Shobuj]{d}{ \Base_n } [swap]{ \circ F \paran{ d_{n+1,0} } } \\ \Hom_{\calC} \paran{ F(n); X }
\end{tikzcd} , \quad 
\forall n\in \num_0 ,
\end{equation}
whose connecting morphisms are provided by right-composition with the zeroeth face maps. This arrangement in \eqref{eqn:def:Base} will be used to formulate a purely categorical notion of convexity. 


At this point, we recall the concept of \textit{decalage} (French for "shift" or "gap") [see \cite[III.5]{GoerssJardine2009} or \cite{Illusie1972complexe}], which is a functorial operation that systematically shifts indices, filtrations, or topological degrees. Its primary utility lies in resolving structural misalignments—either by shifting the pages of a spectral sequence in homological algebra or by providing economical models for path spaces in simplicial homotopy theory. Decalage provides an useful tool for constructions in the theory of \emph{decomposition spaces} \cite[e.g.]{GalvezCarrilloEtal2018a, GalvezCarrilloEtal2018b} and of \textit{bisimplicial sets} \cite[e.g.]{Stevenson2012decalage, Stevenson2016decalage}. The present article uses it for a categorical notion of convexity.

\paragraph{Convexity} Convexity can be given a purely structural interpretation in arbitrary categories, which agrees with the usual notion of convexity in $\Topo$. An object $X$ will be called \emph{convex} if it has a natural transformation $\Cone^X$ which plays the role of a right inverse of $\Base$, as shown in the diagram below : 
\begin{equation} \label{eqn:def:ConeNat}
\begin{tikzcd} [scale cd = 1.5]
	    & && & \Delta^{op} \arrow{drrr}{\Nerve_{F}(X)} \\
\Delta^{op} \arrow{ddr}[swap]{ \Nerve_{F}(X) } & && & {} &&& \SetCat \\
& {} && \tilde{\Delta}^{op} \arrow[dashed, bend right=10, Akashi]{dll}[name = n3]{} \arrow{ulll}[swap]{\iota} \arrow[dashed, Akashi]{urrrr}[name=n1]{} \arrow[dd, pos=0.3, "\text{Shift}"'] \arrow[dashed, Akashi]{drrrr}[name=n2]{} \arrow[uur, pos=0.8, "\iota"] \\ 
& \SetCat && {} & {} &&& \SetCat \\
& && \tilde{\Delta}^{op} \arrow[r, "\iota"] & \Delta^{op} \arrow{urrr}[swap]{\Nerve_{F}(X)}
\arrow[shorten <=1pt, shorten >=1pt, Rightarrow, Shobuj, to path={(n2) to[out=90,in=-45] node[xshift=20pt]{$\Base$} (n1)} ]{  }
\arrow[shorten <=1pt, shorten >=1pt, Rightarrow, Shobuj, to path={(n3) to[out=-90,in=-90] node[yshift=15pt]{$\Cone^X$} (n2)} ]{  }
\arrow[shorten <=1pt, shorten >=1pt, Rightarrow, Shobuj, to path={(n1) to[out=135,in=90] node[yshift=15pt]{$\Id$} (n3)} ]{  }
\end{tikzcd}    
\end{equation}
The thick arrows in \eqref{eqn:def:ConeNat} represent natural transformations. This arrangement is one of the major contributions of this paper. Convexity is usually interpreted as an intrinsic property of a space, and resulting form a closure property with respect to convex linear sums. Such a linear structure may not be available for a general object or category. Equation \eqref{eqn:def:ConeNat} redefines convexity as a relational property. These relations are encoded within the property of the transformation $\Cone$ being \emph{natural}. For the topological example of \eqref{eqn:def:StndrdTopoSmplx}, $\Cone$ becomes the \emph{Cone-construction}. There its connecting morphisms are maps $\Cone^X_n : \Hom \paran{ F(n) ; X } \to \Hom \paran{ F(n+1) ; X }$. Recall that the standard topological $n$-simplex can be described by coordinates 
\[ F_{stndrd, Topo} (n) := \SetDef{ \paran{ t_0, \ldots, t_n } }{ t_j \geq 0, \, \sum_{j=0}^{n} t_j = 1 } , \quad \forall n\in \num_0.\]
Let $a_n$ be an arbitrary point in the $n$-th simplex. Using these coordinates, the action of the (topological) cone transformations may be described as
\begin{equation} \label{eqn:Cone_Topo}
\paran{ \Cone^X_n \sigma } \paran{ t_0, \ldots, t_n, t_{n+1} } :=  
\begin{cases}
	t_{0} a_n + \paran{1-t_{0}} \sigma \paran{ \frac{t_1}{1-t_{0}} , \ldots, \frac{t_{n+1}}{1-t_{0}} } & \mbox{ if } t_{0} < 1 \\
	a_n & \mbox{ if } t_{0}  =1 
\end{cases}, 
\quad \forall n\in \num_.
\end{equation}
The commutation in \eqref{eqn:def:Base} follows from \eqref{eqn:smplc_id}~(i). Thus the $\Base$ transformation exists and is independent of $X$. We next examine the consequences of \eqref{eqn:def:ConeNat}. It implies the existence of a special collection of maps
\begin{equation} \label{eqn:cone:a}
\Cone^X_n : \Hom_{\calC} \paran{ F(n); X } \to \Hom_{\calC} \paran{ F(n+1); X } , \quad n\in \num_0 .
\end{equation}
which satisfy 
\begin{equation} \label{eqn:ConeNat:3}
\begin{tikzcd}
	\Hom_{ \calC } \paran{ F(n); X  } \arrow[d, "\Cone^{X}_n"'] \arrow{rrr}{ \circ F \paran{ d_{n,i} } } &&& \Hom_{ \calC } \paran{ F(n-1); X  } \arrow[d, "\Cone^{X}_{n-1}"] \\ 
	\Hom_{ \calC } \paran{ F(n+1); X  } \arrow{rrr}{ \circ F \paran{ d_{n+1,i+1} } } &&& \Hom_{ \calC } \paran{ F(n); X  }
\end{tikzcd} , \quad 
0\leq i\leq n
\end{equation}
and
\begin{equation} \label{eqn:Base:2}
\begin{tikzcd}
	\Hom_{\calC} \paran{ F(n); X } \arrow{d}[swap]{ \Cone^X(n) } \arrow[r, "="] & \Hom_{\calC} \paran{ F(n); X } \\ 
	\Hom_{\calC} \paran{ F(n+1); X } \arrow[ur, "\circ F\paran{ d_{n+1,0} }"']
\end{tikzcd},
\quad \forall n\in\num_0 .
\end{equation}
Now observe that \eqref{eqn:ConeNat:3} and \eqref{eqn:Base:2} together imply 
\begin{equation} \label{eqn:P:6}
\forall \begin{tikzcd} F(n) \arrow[ "\sigma"', d ] \\ X	\end{tikzcd} , \quad 
\Cone^{X}_n \paran{ \sigma } \circ F\paran{ d_{n+1,i} } \,=\, 
\begin{cases}
	\Cone^{X}_{n-1} \paran{ \sigma \circ F \paran{ d_{n,i-1} } } & \mbox{ if } 1\leq i \leq n+1 \\
	\sigma & \mbox{ if } i=0
\end{cases} ; \quad \forall n\in \num_0 .
\end{equation}
To summarize :

\begin{proposition} [Equivalent characterizations of convexity] \label{prop:cdi9}
The following are equivalent for any object $X$ 
\begin{enumerate} [(i)]
	\item $X$ is convex. 
	\item $X$ has associated to it a cone construction in the sense of \eqref{eqn:def:ConeNat}.  
	\item $X$ has a sequence of maps as in \eqref{eqn:cone:a}, such that the identities in \eqref{eqn:ConeNat:3} and \eqref{eqn:Base:2} are satisfied.
	\item $X$ has a sequence of maps as in \eqref{eqn:cone:a}, such that the identities in \eqref{eqn:P:6} are satisfied. 
\end{enumerate}
\end{proposition}

The original work \cite{Das2024hmlgy} requires the strong assumption :

\begin{Axiom} \label{A:convex}
All the cells $\SetDef{F(n)}{n\in\num_0}$ are convex.
\end{Axiom}

Note that in $\Topo$ the standard topological simplex 
\begin{equation} \label{eqn:def:StndrdTopoSmplx}
F_{stndrd, Topo} : \Delta \to \Topo, \quad F_{stndrd, Topo} (n) = \mbox{ convex span of } \braces{ e^{(n+1)}_0, \ldots, e^{(n+1)}_n } ,
\end{equation}
satisfies the Axiom of convexity. Here $e^{(n+1)}_0, \ldots, e^{(n+1)}_n$ are $n+1$ independent eigenvectors in $\real^{n+1}$. The resulting homology is the singular homology for topological spaces. Note that the $n$-th cell of this functor is isomorphic to the $n$-dimensional closed disk. The need for Axiom \ref{A:convex} emerges from the following result :

\begin{proposition} \label{prop:hmlgy}
\cite[Thm 1]{Das2024hmlgy}  Any functor $F$ as in \eqref{eqn:functorF} creates a homology functor for the category $\calC$. Now suppose that $F$ satisfies the homotopy axioms Axioms \ref{A:C}, \ref{A:brace}, \ref{A:swap}, \ref{A:F1_join}, \ref{A:1_0_cell} and the convexity Axiom \ref{A:convex}. Then homology is homotopy invariant.
\end{proposition}

This completes  a review of the framework built in \cite{Das2024hmlgy} for realizing a homotopy-invariant notion of homology in a general category. Among all the axioms involved, the convexity Axiom \ref{A:convex} appears to be most non-intuitive and difficult to verify.  The next few sections will present a constructive procedure by which this Axiom is satisfied.
\section{The wedge construction} \label{sec:wedge}

\begin{figure}
\centering
\begin{tikzpicture}[scale=0.6, transform shape, framed, background rectangle/.style={double, ultra thick, draw=gray, rounded corners}]
	\input{\Diagrams simplex_recipe.tex}
\end{tikzpicture}
\caption{Outline of the paper. The diagram represents the axiomatic approach of the paper towards building cells in an arbitrary category $\calC$. The diagrams several concepts and constructions, and arrows represent the logical dependencies. The boxes are shaded based on their position in the hierarchy. The starting points are the Axioms in the transparent boxes. The constructions create a sequence of objects called \emph{cells} are analogous to simplexes in Topology. They are convex, contractible and the basis for homology and homotopy for $\calC$. The end result shown in green is the invariance of the resultant homology with respect to homotopy. }
\label{fig:outline:1}
\end{figure} 

The following categorical assumption provides the key principle underlying our construction.

\begin{Axiom} \label{A:pushout}
The category $\calC$ has pushouts.
\end{Axiom}

Recall that for a pair of morphisms with a common domain as shown below on the left,
\[\begin{tikzcd} a \arrow[r, "f"] \arrow[d, "g"'] & b  \\ c  \end{tikzcd}
\begin{tikzcd} {} \arrow[rr, mapsto] && {} \end{tikzcd}
\begin{tikzcd} a \arrow[r, "f"] \arrow[d, "g"'] & b \arrow[d, dotted] \\ c \arrow[r, dotted] & d \end{tikzcd}\]
its pushout is an object $d$ along with morphisms which create the commuting square on the right above.     Moreover, this square is universal in the sense that any other commutation square that extends the wedge diagram on the left, can be factorized uniquely through the commutation square on the right. In $\Topo$ pushouts corresponds exactly to the notion of topological gluing. These various categorical notions - terminal object, products and pushouts, are examples of limits and colimits, known collectively as \emph{universal} constructions in Category theory. They are defined as minimal or final cones and co-cones to certain diagrams.

\paragraph{Wedge construction} Suppose Axioms \ref{A:C}, \ref{A:pushout} and \ref{A:brace} hold. The three ingredients provided by these Axioms - a point-like object $F_0$, an interval-like object $F_1$ and the capability of pushouts - is all that is needed for the constructive procedure. Given any object $X$ of $\calC$, consider the diagram on the left below :
\begin{equation} \label{eqn:def:Wedge}
\begin{tikzcd} [scale cd = 0.8]
	& & 1_{\calC} \\
	X \arrow[r, "\cong"] & X \times 1_{\calC} \arrow[ur, "!"] \arrow[dr, " X \times \tilde{d}_{1,1} "'] \\
	& & X \times F_1
\end{tikzcd} \;
\begin{tikzcd} {} \arrow[rr, mapsto, "\text{pushout}"] && {} \end{tikzcd} \; 
\begin{tikzcd} [scale cd = 0.8]
	& & 1_{\calC} \arrow[dr, Itranga, "\Top_X"] \\
	X \arrow[r, "\cong"] & X \times 1_{\calC} \arrow[ur, "!"] \arrow[dr, " X \times \tilde{d}_{1,1} "'] && \Shobuj{ \Wedge(X) } \\
	& & X \times F_1 \arrow[ur, Itranga, "\Pinch_X"' ]
\end{tikzcd} 
\end{equation}
The diagram on the right above is the resulting pushout. The colimit object $\Wedge(X)$ shown in green will be called the \emph{wedge} construct of $X$. The two additional arrows needed to complete the pushout square have been drawn in red , and labeled $\Pinch_X$ and $\Top_X$ respectively. They will be of paramount importance in uncovering the structure of $\Wedge(X)$. We first argue that this transformation is functorial :

\begin{lemma} \label{lem:Wedge}
If Axioms \ref{A:C}, \ref{A:brace} and \ref{A:pushout} hold, then the pushout construction in \eqref{eqn:def:Wedge} creates a functor $\Wedge : \calC \to \calC$.
\end{lemma}

Lemma \ref{lem:Wedge} is proved in Section \ref{sec:proof:Wedge}. When $\calC=\Topo$, $\Wedge(X)$ is the space obtained by first creating a cylinder $X\times I$ and then pinching the top face into a single point. This topological interpretation is the motivation behind naming the morphism $\Pinch$ in \eqref{eqn:def:Wedge}. The pinch morphism connects the cylinder to the cone. The original object $X$ attaches via the boundary maps $\tilde{d}_{1,0}, \tilde{d}_{1,1}$ to both the bottom and top of the cylinder. The former attachment can be used in composition with $\Pinch$ to create a concept of a \emph{floor} or \emph{bottom-face} for the wedge product : 
\begin{equation} \label{eqn:def:Bottom}
\begin{tikzcd}
	X \times F_0 \arrow[rr, " X \times \tilde{d}_{1,0} "] & & X \times F_1 \arrow[d, "\Pinch_X"] \\
	X \arrow[u, "="] \arrow[rr, Shobuj, " \Floor_X"', dashed] & & \Wedge(X)
\end{tikzcd}
\end{equation}
The next observation is

\begin{lemma} \label{lem:kd0d3}
If Axioms \ref{A:C}, \ref{A:brace} and \ref{A:pushout} hold, then for every contractible object $X\in \calC$, there is a morphism $\Flatten_X : \Wedge(X) \to X$ which is also the left inverse of the Bottom-morphism \eqref{eqn:def:Bottom}.
\end{lemma}

Lemma \ref{lem:kd0d3} essentially says that if $X$ is contractible then $\Wedge(X)$ can be retracted onto $X$. Note that the condition of contractibility cannot be dropped from the statement. In fact the Bower fixed point theorem states that this is not true when $\calC = \Topo$ and $X = S^1$, a contractible space. In that case $\Wedge(X) = D^2$, the 2-disk. Lemma \ref{lem:kd0d3} can be proved by construction :
\begin{equation} \label{eqn:conractible:1}
\begin{tikzcd} [scale cd = 1.5]
	X \arrow[ddrrrr, bend right = 30, "\Floor_X"' ] \arrow[rr, "{X \times \tilde{d}_{1,0}}"] \arrow[rrrrd, "="'{pos=0.5}, bend right] & & X\times F_1 \arrow[rrd, "\Contraction_X"{pos=0.3}, bend right, Holud] \arrow[rrdd, "\Pinch_X"'{pos=0.1}, bend right=49, Akashi] & & X\times F_0 \arrow[r, "!"] \arrow[ll, "{X \times \tilde{d}_{1,1}}"'] & F_0 \arrow[ld, "x", Holud] \arrow[ldd, pos=0.3, "\Top_X", bend left=49, Akashi] \\
	& & & & X & \\
	& & & & \Wedge(X) \arrow[u, "\Flatten_X"', Shobuj] & 
\end{tikzcd}
\end{equation}
The commuting square created by the two yellow loops are lifted from the definition of contractibility in \eqref{eqn:conractible:3}. The yellow arrows can be seen to create a cocone for the diagram on the left of \eqref{eqn:def:Wedge}. The wedge construct is by deifnition the limiting cocone for this diagram. Thus the yellow cocone must be factorizable via the wedge cocone. The connecting morphism that makes this possible is drawn in green and labelled $\Flatten_X$. The left most loop created by $\Floor_X$ comes from its definition \eqref{eqn:def:Bottom}. The final loop loop to verify is the top left, created by $\Contraction_X$ and $\Id_X$. This is again lifted directly from \eqref{eqn:conractible:3}. Thus \eqref{eqn:conractible:1} is a consistent commutative diagrams, showing the various inter-relations between the various special morphisms we have created so far. The bottom left commutation loop which is automatically created, can be written as :
\begin{equation} \label{eqn:conractible:2}
X \mbox{ contractible } \imply \Flatten_X \circ \Floor_X = \Id_X ,
\end{equation}
which is precisely the second claim of Lemma \ref{lem:kd0d3}. When $\calC = \Topo$, $\Flatten(X,x)$ squashes or flattens the wedge construct $\Wedge(X)$ onto its bottom base $X$ in such a way that the crest of the wedge is mapped into the point $x\in X$.

The following axiom on the beginning object of our construction will be instrumental in establishing contractibility throughout the rest of the procedure :

\begin{Axiom} \label{A:1_contract}
The $1$-cell $F_1$ is contractible.
\end{Axiom}

Axiom \ref{A:1_contract} is essential for the following lemma, which forms the underlying principle of our construct.

\begin{lemma} \label{lem:f3sk0}
If Axioms \ref{A:C}, \ref{A:brace}, \ref{A:F1_join} and \ref{A:1_contract} hold, then for any object $X$, its wedge construct $\Wedge(X)$ is contractible.
\end{lemma}

Lemma \ref{lem:f3sk0} is proved in Section \ref{sec:proof:f3sk0}. Once again, when $\calC = \Topo$ one can verify that the mapping cone is not guaranteed to be contractible unless the cylinder created uses the a contractible object like the interval. Lemma \ref{lem:f3sk0} provides sufficient conditions under which an object created as a wedge is contractible. Eventually, our cells will be created by iterations of the wedge construct. Lemma \ref{lem:f3sk0} will then be instrumental is concluding that the cells are contractible. Being contractible leads to the existence of the morphisms $\Flatten$ and $\Contraction$ from Diagram \eqref{eqn:conractible:1}. These morphisms will play the role of faces and degeneracies in the construction of our co-simplicial object $F$.

The next section presents how these few and simple collection of morphisms generate the full set of face and degeneracy maps.

\section{Faces of the simplex} \label{sec:induction}

Recall that our task is to define a functor $F$ on the morphisms in $\Delta$ satisfying the various assumptions made in Section \ref{sen:hmlgy_hmtpy}.  In general, it is not easy to make an explicit construction of a functor of the form \eqref{eqn:functorF}, as the morphisms incoming at $n$ or outgoing at $n$ increase exponentially with $n$. The task becomes easier if one concentrates only on a special collection of morphisms known as \emph{face maps} :
\begin{equation} \label{eqn:def:face}
0\leq i\leq n \;:\; \Shobuj{d_{n,i}} : [n-1] \to [n], \, \quad j \mapsto \begin{cases}
	j & \mbox{ if } j<i	\\
	j+1 & \mbox{ if } j\geq i	
\end{cases}.
\end{equation}
and \emph{degeneracy} maps :
\begin{equation} \label{eqn:def:degeneracy}
0\leq i\leq n \;:\;  \Shobuj{s_{n,i}} : [n+1] \to [n], \quad j \mapsto \begin{cases}
	j & \mbox{ if } j \leq i	\\
	j-1 & \mbox{ if } j > i	
\end{cases}.
\end{equation}
To put more concisely, the $i$-th face map $d_{n,i}$ is the unique map which skips the element $i$ in $[n+1]$, and the $i$-th degeneracy map $s_{n,i}$ is the unique surjective map such that the pre-image of element $i$ has two elements. The face and degeneracy maps generate every morphism $\phi:[m]\to [n]$ in $\Delta$. The obey the following set of identities :
\begin{equation} \label{eqn:smplc_id}
\begin{split}
	d_{n+1,i} d_{n,j} = d_{n+1,j+1} d_{n,i} & \mbox{ if } i<j \\
	s_{n-1,j} d_{n,i} = d_{n-1,i} s_{n-2,j-1} & \mbox{ if } i<j \\
	s_{n-1,j} d_{n,j} = \Id_{[n-1]} & \\
	s_{n-1,j} d_{n,j+1} = \Id_{[n-1]} \\
	s_{n-1,j} d_{n,i} = d_{n-1,i-1} s_{n-2,j} & \mbox{ if } i>j+1 \\
	s_{n-1,j} s_{n,i} = s_{n-1,i} s_{n,j+1} & \mbox{ if } i\leq j 
\end{split} 
\end{equation}
which are called the \emph{simplicial identities}. Equations \eqref{eqn:def:face}, \eqref{eqn:def:degeneracy} and \eqref{eqn:smplc_id} provide a sufficient set of morphisms and composition relations that generate all the morphisms of $\Delta$. This means that given a collection of morphisms \eqref{eqn:def:face} and \eqref{eqn:def:degeneracy} which satisfy \eqref{eqn:smplc_id}, any morphism in $\Delta$ may be expressed as a composition of different $d_{n,i}$ and $s_{n',i'}$-s. Thus when trying to construct a functor $F$ as in \eqref{eqn:functorF}, one only needs to specify the objects $\SetDef{ F(n) }{ n\in \num_0 }$, the action of $F$ on the face and degeneracy morphisms, and ensure that $F$ preserves the simplicial identities of \eqref{eqn:smplc_id}. 

\paragraph{Cells} We first define the action of $F$ on $\num_0$, i.e., specify the cells. The cells are created iteratively via the wedge construction : 
\begin{equation} \label{eqn:Wedge_cells:1}
\forall n\in\num  \,: \quad F_{n+1} := \Wedge \paran{ F_n } .
\end{equation}

\paragraph{Face maps} Recall that $F_0$ denotes the terminal object $1$. To create face maps we redraw \eqref{eqn:def:Wedge} but for the cells : 
\begin{equation} \label{eqn:Wedge_cells}
\begin{tikzcd} [scale cd = 0.8]
	& & F_0 \\
	F_n \arrow[r, "\cong"] & F_n \times F_0 \arrow[ur, "!"] \arrow[dr, " F_n \times \tilde{d}_{1,1} "'] \\
	& & F_n \times F_1
\end{tikzcd} \;
\begin{tikzcd} {} \arrow[rr, mapsto, "\text{pushout}"] && {} \end{tikzcd} \; 
\begin{tikzcd} [scale cd = 0.8]
	& & F_0 \arrow[dr, Itranga] \\
	F_n \arrow[r, "\cong"] & F_n \times F_0 \arrow[ur, "!"] \arrow[dr, " F_n \times \tilde{d}_{1,1} "] && F_{n+1} \\
	& & F_n \times F_1 \arrow[ur, Itranga, "\Pinch_{F_n}"'] 
\end{tikzcd} 
\end{equation}
The $\Floor$ morphism helps us define the face map $\tilde{d}_{n+1,n+1}$  : 
\begin{equation} \label{eqn:Wedge_cells:2}
\forall n\in\num  \,: 
\begin{tikzcd} F_{n} \arrow[d, " \tilde{d}_{n+1,j} "', Shobuj] \\ F_{n+1} \end{tikzcd} := 
\begin{cases}
	\Wedge \paran{ \tilde{d}_{n,j} } & \mbox{ if } 0\leq j< n+1 \\
	\Floor_{F_n} \mbox{ \eqref{eqn:def:Bottom}} & \mbox{ if } j=n+1
\end{cases}
\end{equation}

\paragraph{Centroids} Henceforth let us assume that the cells $F_n$ are contractible. Then we have morphism $\Flatten(X,x)$ from \eqref{eqn:conractible:1} at our disposal. They will be used to construct some of the degeneracy maps. For that purpose we need elements $a_n$ in each cell $F_N$. This will be guaranteed by Axiom \ref{A:F1_join}. The diagram on the left below is a repeat of the pushout diagram representing the colimit of Axiom \ref{A:F1_join} :
\begin{equation} \label{eqn:def:a1}
\begin{tikzcd}
	& F_1 \arrow[dr, Itranga] \\
	1_{\calC} \arrow[rr, Shobuj, "a_1", dashed] \arrow[ur, "\tilde{d}_{1,0}"] \arrow[dr, "\tilde{d}_{1,1}"'] && F_1 \\
	& F_1 \arrow[ur, Itranga] \\
\end{tikzcd} , \quad 
\begin{tikzcd}
	& F_1 \\
	1_{\calC} \arrow[dashed, Shobuj, dd] \arrow[dr, "a_n"'] \arrow[ur, "a_1"] \arrow{rr}{ \bracket{ a_1, a_n } } && F_{n} \times F_{1} \arrow[rrr, "\Contraction_{F_n}"] \arrow[dl, "\proj_1"] \arrow[ul, "\proj_2"'] &&& F_{n+1} \\ 
	& F_n \\
	{} \arrow[dashed, Shobuj, "a_{n+1}", rrrrr] &&&&& {} \arrow[dashed, Shobuj, uu]
\end{tikzcd}
\end{equation}
The morphism $\text{center}$ in the statement of Axiom \ref{A:F1_join} has been re-written as a well defined element $a_1$ of $F_1$. The diagram on the right shows how the elements $a_n$ of cells $F_n$ are created inductively. The point $a_{n+1}$ is the image of the point $\bracket{ a_1, a_n }$ in the product space, under the map $\Pinch_{F_n}$ introduced in \eqref{eqn:Wedge_cells}.

So far we have a sequence of cells $F_n$ defined inductively in \eqref{eqn:Wedge_cells:1}, and points $a_n \in F_n$ defined inductively in \eqref{eqn:def:a1}. We can now define the degeneracy maps as 		
\begin{equation} \label{eqn:Wedge_cells:3}
\forall n\in\num  \,: 
\begin{tikzcd} F_{n+1} \arrow[d, " \tilde{s}_{n,j} "', Shobuj] \\ F_{n} \end{tikzcd} := 
\begin{cases}
	\Wedge \paran{ \tilde{s}_{n-1,j} } & \mbox{ if } 0\leq j< n \\
	\Flatten \paran{ F_n, a_{n} } \mbox{ Equation \eqref{eqn:conractible:1}} & \mbox{ if } j=n
\end{cases}
\end{equation}
We have this seen an iterative procedure by which a sequence of objects $F_0, F_1, F_2, \ldots$ have been created, along with morphisms $\tilde{d}_{n,j} : F_{n-1} \to F_n$ and $\tilde{s}_{n,j} : F_{n+1} \to F_{n}$. We state our main result in the next section, where we verify that these objects $F_n$ are proxies for the images $F(n)$ of a co-simplicial object as in \eqref{eqn:functorF}, and the $\tilde{s}_{n,j}$ and $\tilde{d}_{n,j}$ are proxies for the boundary and degeneracy morphisms $s_{n,j}$ and $d_{n,j}$ respectively.


We are now ready to describe how the various objects and morphisms assemble to construct the functor \eqref{eqn:functorF} from Theorem \ref{thm:recipe}~(i).

\begin{theorem} \label{thm:recipe}
Suppose Axioms Axioms \ref{A:C}, \ref{A:swap}, \ref{A:F1_join}, \ref{A:pushout}, \ref{A:brace}, and \ref{A:1_contract} hold.
\begin{enumerate} [(i)]
	\item The objects $F_n$ defined inductively in \eqref{eqn:Wedge_cells:1}, morphisms $\tilde{d}_{n+1,j}$ defined in \eqref{eqn:Wedge_cells:2}, and morphisms $\tilde{s}_{n+1,j}$ defined in \eqref{eqn:Wedge_cells:3} assemble to create a functor as in \eqref{eqn:functorF}, which satisfies Axiom \ref{A:1_0_cell}. 	
	\item All the cells $F_n$ are contractible.		
	\item If $F_1$ is convex, then Axiom \ref{A:convex} is satisfied. 
\end{enumerate}
\end{theorem}

Thus the six categorical axioms -- Axioms \ref{A:C}, \ref{A:swap}, \ref{A:F1_join}, \ref{A:pushout}, \ref{A:brace}, \ref{A:1_contract}, create a generation scheme for cells, faces and degeneracies. The ingredients of the generative procedure are simply the 0 and 1 cells $F_0 = 1_{\calC}$ and $F_1$ respectively. As promised by claim (ii) of the theorem, all the cells are also convex. The rules for attachment of the cells are provided by the trivial endpoint morphisms of Axiom \ref{A:brace} as well as the pushout Axiom \ref{A:F1_join}. 

The rest of this section contains the proof of Theorem \ref{thm:recipe}.

\paragraph{Proof of contractibility} Since $F_1$ is contractible, each of the cells are contractible as a direct consequence of Lemma \ref{lem:f3sk0} and the recursive nature of the definition. This proves the claim (ii) of Theorem \ref{thm:recipe}. Before proving functoriality, we establish convexity.

\paragraph{Cone construction for cells} The key ingredient of Theorem \ref{thm:recipe} is the creation of a cone construction. We thus have to create a family of maps as in \eqref{eqn:cone:a}. In our case they take the form :
\[ \Cone^{F_m}_n : \Hom_{ \calC } \paran{ F_n; F_m } \to \Hom_{ \calC } \paran{ F_{n+1}; F_m } , \quad m,n\in \num_0 ,   \]
The construction of $\Cone^{F_m}_n$ is via 	
\begin{equation} \label{eqn:def:Cone_recipe}
\begin{tikzcd} F_n \arrow[d, "\sigma"'] \\ F_m \end{tikzcd} \;\Rightarrow\; 
\begin{tikzcd} \Wedge(F_n) \arrow[d, "\Wedge( \sigma )"'] \\ \Wedge(F_m) \end{tikzcd} \;\Rightarrow\; 
\begin{tikzcd} [column sep = large]
	F_{n+1} \arrow[Shobuj, dashed, d, "\Cone^{F_m}_n (\sigma) "] \arrow[r, "="] & \Wedge(F_n) \arrow[d, "\Wedge( \sigma )"] \\
	F_m & \Wedge(F_m) \arrow[l, " \tilde{s}_{m,m} "]
\end{tikzcd}
\end{equation}
By taking advantage of Proposition \ref{prop:cdi9}, it is enough to show that the identities in \eqref{eqn:P:6} are satisfied. If $i<n+1$,
\[\begin{split}
\Cone^{F_m}_{n} \paran{ \sigma } \circ \tilde{d}_{n+1,i} &= \tilde{s}_{m,m} \circ \Wedge(\sigma) \circ \tilde{d}_{n+1,i}, \quad \mbox{ by \eqref{eqn:def:Cone_recipe} } \\
&= \tilde{s}_{m,m} \circ \Wedge(\sigma) \circ \Wedge \paran{ \tilde{d}_{n,i} }, \quad \mbox{ by \eqref{eqn:Wedge_cells:2} } \\
&= \tilde{s}_{m,m} \circ \Wedge \paran{ \sigma \circ \tilde{d}_{n,i} } \\
& = \Cone^{F_m}_{n-1} \paran{ \sigma \circ \tilde{d}_{n,i} }, \quad \mbox{ by \eqref{eqn:def:Cone_recipe} } .
\end{split}\]
This verifies the first identity of \eqref{eqn:P:6}. When $i=n+1$,
\[\begin{split}
\Cone^{F_m}_{n} \paran{ \sigma } \circ \tilde{d}_{n+1,n+1} &= \tilde{s}_{m,m} \circ \Wedge(\sigma) \circ \tilde{d}_{n+1,n+1}, \quad \mbox{ by \eqref{eqn:def:Cone_recipe} } \\
&= \Flatten \paran{ F_m, a_{m} } \circ \Wedge(\sigma) \circ \Floor_{F_n}, \quad \mbox{ by \eqref{eqn:Wedge_cells:2}, \eqref{eqn:Wedge_cells:3} }
\end{split}\]
This completes the verification of \eqref{eqn:P:6} for the cone constructed as in \eqref{eqn:def:Cone_recipe}. Thus the cells are indeed convex, as claimed in Theorem \ref{thm:recipe}~(iii). It remains to prove Claim (i).

\paragraph{Proof of functoriality} The face and degeneracy morphisms have been defined recursively in \eqref{eqn:Wedge_cells:2} and \eqref{eqn:Wedge_cells:3} using the Wedge construction. The Wedge construction has been shown to be functorial in Lemma \ref{lem:Wedge} and therefore preserves composition. Therefore if the simplicial identities \eqref{eqn:smplc_id} hold for some $n$, they would also hold for $n+1$ for all values of the indices $i,j$ except the extremal ones. 

\paragraph{Base case} When we check the base case for the simplicial identities \eqref{eqn:smplc_id}, all the identities with a face map as the left-most term become morphisms to $F_0=1$. Because the identities claim equalities of morphisms between identical pair of objects and $1$ is a terminal object, the base case is automatically satisfied for these identities. The only identity from \eqref{eqn:smplc_id} whose base case needs to be verified is 
\[\tilde{d}_{2,0} \tilde{d}_{1,1} = \tilde{d}_{2,2} \tilde{d}_{1,0} .\]
Both sides of the equality are arrows starting from $1$ and its verification is trivial.

\paragraph{Some identities} Each of the six equations of \eqref{eqn:smplc_id} can be divided into sub-cases depending on the values of the indices $i$ and $j$. Figures \ref{fig:smplc_id_1}, \ref{fig:smplc_id_2}, \ref{fig:smplc_id_3}, \ref{fig:smplc_id_4}, \ref{fig:smplc_id_5}, and \ref{fig:smplc_id_6} display the breakdown of each of these equations, and the Lemma or identity that proves each of the sub-cases. The rest of the section proves the identities referred to in these figures. We next present some identities that prove the sub-cases present in these diagrams.
\begin{equation} \label{eqn:1_3_c}
\Wedge\paran{ \Floor_{F_{n-1}} } \Floor_{F_{n-1}} = \Floor_{F_{n}} \Floor_{F_{n-1}} . 
\end{equation}
\begin{equation} \label{eqn:1_3_a}
\Wedge\paran{ d_{n,i} } \Floor_{F_{n-1}} = \Floor_{F_{n}} \Wedge\paran{ d_{n-1,i} }, \quad 0 \leq i < n.
\end{equation}
\begin{equation} \label{2_3_b}
\Flatten\paran{ F_{n-1}, a_{n-1} } \Wedge\paran{ \Floor_{F_{n-2}} } = \Floor_{F_{n-2}} \Flatten\paran{ F_{n-2}, a_{n-2} } ;
\end{equation}
\begin{equation} \label{2_2_b}
\Flatten\paran{ F_{n-1}, a_{n-1} } d_{n,i} = d_{n-1,i} \Flatten\paran{ F_{n-2}, a_{n-2} } , \quad 0 \leq i < n.
\end{equation}
\begin{equation} \label{eqn:5_2_b}
s_{n-1,j} \Floor_{F_{n-1}}  = \Floor_{F_{n-2}} s_{n-2,j} , \quad 0 \leq j< n-2 .
\end{equation}

These identities are sufficient to prove the simplicial identities, as outlined in Figures \ref{fig:smplc_id_1}, \ref{fig:smplc_id_2}, \ref{fig:smplc_id_3}, \ref{fig:smplc_id_4}, \ref{fig:smplc_id_5}, and \ref{fig:smplc_id_6}. These identities are consequences of the wedge construction and are proved in Section \ref{sec:identities}. This completes the proof of Theorem \ref{thm:recipe}~(i). \qed 

This completes the presentation the our main theoretical results. Some instances of our axiomatic framework are presented next.
\section{Examples of the axiomatic framework} \label{sec:example}

As outlined in Figure \ref{fig:outline:4}, a category $\calC$ which satisfies Axioms \ref{A:C}, \ref{A:brace}, \ref{A:swap}, \ref{A:F1_join}, \ref{A:pushout} and \ref{A:1_contract} has a simplicial object $F:\Delta \to \calC$ whose cells are contractible and convex, and thus leads to a homology and homotopy which are mutually compatible. We now examine the some useful realizations of this axiomatic framework, and some examples where these are partially fulfilled.

As displayed in Figure \ref{fig:outline:4}, the ground Axioms are either about the category $\calC$, or about the object $F_1$. We shall call an object $F_1$ to be a \textit{generalized interval} if it satisfies Axioms \ref{A:brace}, \ref{A:F1_join} and \ref{A:1_contract}. It will be called symmetric or asymmetric based on whether it satisfies Axiom \ref{A:swap}.

\paragraph{I. Topology} The prime motivation of this work is of course the category $\Topo$ of topological spaces and continuous maps. The compact interval $I$ plays the role of $F_1$ and the resulting functor $F$ is the standard simplicial space in which the $n$-th cell is the $n$-dimensional simplex. 

A key constraint of our axiomatic framework is Axiom \ref{A:F1_join}. It provides the structural basis for composition of homotopies. In fact the nature of the interval $I$ as a gluing of two copies of itself leads to a categorical notion of \emph{self-similarity} \cite{Leinster2011selfsim, Freyd2008algebraic}. This Axiom remains one of the hardest to be satisfied in other categories.

\paragraph{II. Sets} Since $\Topo$ embeds into $\SetCat$, the category of sets up to a certain cardinality and set-maps, $\SetCat$ is also a trivial instance where our axiomatic framework is satisfied.

\paragraph{III. Topological concrete category} A convenient way to leverage the tools of topology is using the premises of a topological concrete category \cite[e.g.]{Herrlich1974topo1, herrlich1974topo2, brummer1984topological, DikranjanEtAl1988topo}. Such a category extracts certain categorical / structural properties of $\Topo$ to enable a more general and abstract category in which several familiar notions of topology can be realized along with their usual properties. A topological concrete category $\calC$ is equipped with a faithful functor $U: \calC \to \SetCat$ satisfying the following two conditions:
\begin{enumerate} [(i)]
\item For any set $X$ and any family of morphisms in $\SetCat$ indexed by $i \in I$:
\[f_i: X \to U(X_i), \; i\in I,\]
there exists a unique object $X^*$ in $\calC$ such that $U(X^*) = X$, and a family of $\calC$-morphisms $\bar{f}_i: X^* \to X_i$ such that $U(\bar{f}_i) = f_i$.
\item This $X^*$ must be initial, meaning for any object $Y\in \calC$, a map $g: U(Y) \to X$, and morphisms $f'_i : Y \to X_i$ such that $f_i \circ g = U(f'_i)$, then there must exist a $\calC$-morphism $\phi : Y \to X^*$ such that $U(\phi) = g$ and $f_i \circ \phi = f'_i$
\end{enumerate}

\begin{figure} [!t]
\centering
\begin{tikzpicture}[scale=1, transform shape]
	\input{\Diagrams concrete_classification.tex}
\end{tikzpicture}
\caption{A broad classification of categories with generalized interval objects.  Note that the pointed versions of all these categories also have interval, as do functor categories whose codomains are any of these categories. See the discussion in Section \ref{sec:example}.}
\label{fig:cncrt_class}
\end{figure}

These conditions are the categorical generalization of the initial topology (like the subspace or product topology). It says one can always "pull back" structure from a collection of maps $f_i$ into unstructured sets $U(X_i)$. The forgetful functor $U: \calC \to \SetCat$ has the useful property :

\begin{lemma} \label{lem:jd9ps}
\cite[Prop 21.15]{AdamekEtAl2004cats} If $(\calC, U)$ is topological concrete, then $U$ uniquely lifts both limits (via initiality) and colimits (via finality) in $\SetCat$, and it preserves both limits and colimits.
\end{lemma}

The colimit preserving property arises from the fact that $U$ possesses a left adjoint, called the \emph{discretization functor}. Thus if we take the lift of the diagrams 
\[\begin{tikzcd}
\star \arrow[dr, "\text{left}"'] & & \star \arrow[dl, "\text{right}"] \\ & I
\end{tikzcd} ; \quad
\begin{tikzcd}
& I \arrow[dr, dotted] \\
1 \arrow[ur] \arrow[dr] && I \\
& I \arrow[ur, dotted]
\end{tikzcd}\]
in $\SetCat$ into $\calC$, one would get an interval like object in $\calC$ as well. 

While $\Topo$ remains an example of a topological concrete category, other categories such as $\MeasCat$, category of \textit{uniform spaces}, category of \emph{proximity spaces}, category of \emph{convergence spaces}, and Grothendieck toposes are also examples. See Figure \ref{fig:cncrt_class} for a simple classification of categories which allow generalized interval objects. All of these are valid candidates for our framework. 


\paragraph{IV. Posets} The category $\Poset$ comprises of all partially ordered sets, with order preserving maps being the morphisms. This category is equipped with all limits and colimits, and thus has a terminal object, products and pushouts as well. While there is no unique choice of $F_1$, $F_1$ could be any \textit{bounded poset}, i.e. a poset with a maximum and minimum. However, irrespective of the choice of $F_1$ Axiom \ref{A:swap} would be violated as the ordering cannot be reversed. As a result, while one would still have a wedge construction, cells, homology and homotopy, the notion of homotopy would be transitive and reflexive but not symmetric.

\paragraph{V. Directed spaces} The category $\DirTopo$ of directed spaces provides a framework for "directed algebraic topology" \cite[e.g.]{Raussen2010trace, LeeYetter2022strat}, and is another instance of a topological concrete category. The terminal object is the singleton directed space $\{*\}$. One could take for $F_1$ the directed unit interval $(I, dI)$. Here $I$ is the standard topological interval $[0, 1]$, while the d-paths ($dI$) are the continuous non-decreasing maps $\alpha: [0, 1] \to [0, 1]$.

While familiar categories such as $\Topo$ and $\MeasCat$ fit well into our framework, there are some other important categories for which there is no easy realization of these framework, such as the category $\StochCat$ of probability spaces and Markov kernels, or the category of dynamical systems \cite[e.g.]{DasSuda2024recon, DasSuda2025enrich}. These extensions offer an interesting topic of investigation.

\section{Simplicial complexes} \label{sec:complexes}

\begin{SCfigure}[][!t]
\centering
\includegraphics[width=0.6\linewidth]{\Diagrams 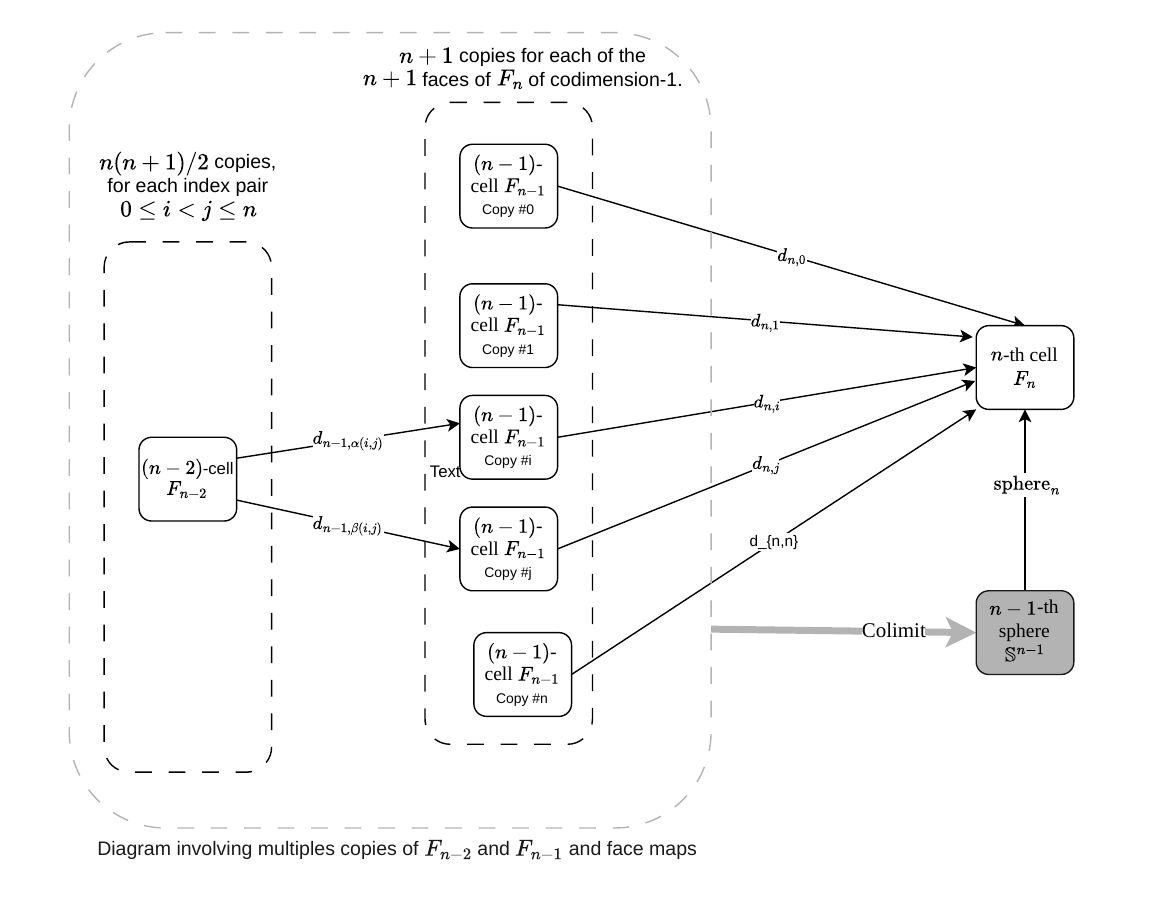}
\caption{Construction of a categorical sphere. The indices $\alpha, \beta$ are defined in \eqref{eqn:sphere_construct}. Each of the boxes correspond to objects in $\calC$ and arrows correspond to morphisms. The layout within the dashed box is a finite diagram in $\calC$. It contains multiple copies of the $(n-1)$ and $n$-cells constructed in Section \ref{sec:induction}, connected with various choices of the boundary morphisms defined in \eqref{eqn:Wedge_cells:2}. The object $\Sphere{n-1}$ is the colimit of this diagram in $\calC$. If $\calC=\Topo$, then $\Sphere{n-1}$ is the usual topological sphere created by viewing the boundary of the $n+1$-simplex $F_n$ as a union of $(n+1)$ pieces of $F_{n-1}$, and gluing these boundary pieces along their common edges $F_{n-1}$. By a reuse of notation, $\SphereMap{n}$ is used to denote the inclusion of $\Sphere{n-1}$ into $F_{n}$. }
\label{fig:sphere_construct}
\end{SCfigure}

Our constructions and theorems performed a "topologification" of the category $\calC$. Topology allows multiple ways of building complicated objects from a countable collection simple building blocks. We can do the same for $\calC$, by using the cells $\braces{F_n}_{n=0}^{\infty}$ for topological simplexes, and pushouts for gluing. We begin with a brief summary of the many ways one builds compound objects in $\Topo$.

\paragraph{Cell complexes}  A cell complex $X$ is a topological space built up inductively by gluing "cells" together. One starts with a discrete set of points ($0$-cells), and inductively, the $n$-th cells are attached along their boundaries to the existing $(n-1)$-skeleton. The $(n-1)$-skeleton is inductively defined to be the union of all cells of dimension up to $(n-1)$. In $\Topo$, these cells are the $n$ simplices isomorphic to the $n$-disks $D^n$, and their boundaries are the spheres $S^{n-1}$. 

The notion of such cell-complexes can be recreated in our abstract framework as well. We begin with a generalization of the notion of a \emph{sphere}. In $\Topo$ a $k$-dimensional sphere is the boundary of $(k+1)$-cell $D^{k+1}$. In a general category, the general cell $F_{k+1}$ takes the place of $D^{k+1}$, and the boundary is a composite of the various face maps. Define the following index set
\begin{equation} \label{eqn:sphere_construct}
\begin{array}{c} \mbox{For each} \\0\leq i < j \leq n \end{array} ,\,
\begin{array}{c} \exists \mbox{ indices} \\0\leq \alpha(i,j), \beta(i,j) <n \end{array}
\mbox{ s.t. }
\begin{array}{c} d_{n,i}d_{n-1,\alpha} \\ = \\ d_{n,j}d_{n-1,\beta}     \end{array}
\end{equation}
Then the sphere $\mathbb{S}^{n-1}$ is a colimit of a diagram as shown in Figure \ref{fig:sphere_construct}. In short, a sphere residing as a boundary of a cell is denoted as the sub-object / monomorphism
\begin{equation} \label{eqn:def:sphere}
\begin{tikzcd} \Sphere{n-1} \arrow[rrr, "\SphereMap{n}", hook] &&& F_n \end{tikzcd} , \quad n=0,1,2,\ldots
\end{equation}
The notion of boundary of a cell enables the important notion of \textit{handle-attachment} or cell-gluing, Given a $\calC$-object $X$ and and integer $k\in \num$, a \emph{$k$-handle attachment} to $X$ is a colimit / pushout of the following diagram
\[\begin{tikzcd}
&& F_k \arrow[drr] \\
\Sphere{k-1} \arrow[urr, "\SphereMap{k}", "\subset"'] \arrow[drr, "\alpha"'] && && X \sqcup_{\alpha} F_k \\
&& X \arrow[urr]
\end{tikzcd}\]
In the diagram above, $\alpha$ is any morphism between the corresponding objects. The resulting colimit object $X \sqcup_{\alpha} F_k$ thus depends on the choice of $\alpha$. One can attach an arbitrary number of $k$-handles in this manner to $X$. The construction can be expressed by a similar colimit diagram
\[\begin{tikzcd}
& & && F_k \arrow[dddrr] \\
& & && \vdots \\
& & && F_k \arrow[drr] \\
\Sphere{k-1} \arrow[uuurrrr, "\SphereMap{k}", "\subset"'] \arrow[drrrr, "\alpha_m"'] & \cdots & \Sphere{k-1} \arrow[urr, "\SphereMap{k}", "\subset"'] \arrow[drr, "\alpha_1"] && && X \sqcup_{\alpha_1} F_k \cdots \sqcup_{\alpha_m} F_k \\
& & && X \arrow[urr]
\end{tikzcd}\]
A \textit{CW complex} or \emph{cellular complex} \cite{Whitehead1949combin} in $\calC$ is constructed by a sequential construction represented as 
\[ X_{0} \xrightarrow{\subset} X_{1} \xrightarrow{\subset} \cdots X_{k-1} \xrightarrow{\subset} X_{k} \xrightarrow{\subset} \cdots \]
The object $X_k$ is created recursively out of $X_{k-1}$ by taking a collection of $k$-cells $\SetDef{ F^{(i)}_k }{i\in I}$ and gluing each of them to $X_{k-1}$. These maps are called \emph{attachment maps}. The associated cellular complex $X$ is the colimit of the above sequence. If the sequence is finite and ends at some $X_N$, $N$ is called the dimension of the cellular object. The $k$-th object $X_k$ is called the $k$-th skeleton.

\paragraph{CW complexes} A CW complex in $\Topo$ is essentially a well-behaved cell complex. The distinction only appears when the space $X$ is built from an infinite number of cells. Formally, it is a cell-complex satisfying two additional conditions, introduced by Whitehead \cite{Whitehead1949combin, Whitehead1949combin2} to rule out pathological topologies.
\begin{enumerate} [(i)]
\item (C) Closure-finiteness : The boundary of every cell is covered by only a finite number of lower-dimensional cells. Thus one cannot have a single cell whose boundary stretches across infinitely many other cells. 
\item (W) Weak topology : A subset of the total space is declared to be closed (or open) if and only if its intersection with the closure of every single cell is closed (or open) within that cell. This forces the global topology of the space to be completely determined by the local topology of its individual cells. 
\end{enumerate}
The name "CW" stands for these two conditions. In a space with a finite number of cells (like a sphere, a torus, or a finite graph in), closure-finiteness and Weak topology are automatically satisfied. In this case, a cell complex and a CW complex are the exact same thing. There is no neat and direct way to accommodate Axiom (W) into a general category $\calC$ such as ours.

\paragraph{Geometric simplicial complexes} These are cell complexes in which the building blocks are not just arbitrary topological disks as they are in a standard cell complex, but rigid, geometric simplices—points ($0$-simplices), line segments ($1$-simplices), solid triangles ($2$-simplices), solid tetrahedra ($3$-simplices), and so on. Because the gluing maps are strict identity maps along the boundaries, a geometric simplicial complex (G.S.C.) ends up being bound by just two rules  - 
\begin{enumerate} [(i)]
\item Every face of a simplex in the complex is also in the complex.
\item The non-empty intersection of any two simplices is exactly one common face.
\end{enumerate}
G.S.C.s are purely combinatorial in nature due to a lack of choice in the attachment maps. They are very amenable to our categorical construction. A finite category $J$ will be called a \emph{cell-complex diagram} if it satisfies the following properties : 
\begin{enumerate} [(i)]
\item The objects of $J$ can be partitioned into levels $J^{(1)}, \ldots, J^{(L)}$.
\item The morphisms in $J$ are generated by arrows that connect an $l$-level object to an $(l+1)$-level object.
\end{enumerate}
A \emph{ G.S.C. } in $\calC$ is a functor from $\Phi : J \to \calC$ such that
\begin{enumerate} [(i)]
\item $\Phi$ maps each $l$-level object into the object $F(l)$.
\item Each morphism $f$ in $J$ that connects an $l$-level object to an $(l+1)$-level object, is mapped into a face map $F \paran{ d_{l+1,j} }$ for some $0\leq j \leq l+1$.
\end{enumerate}
Thus a G.S.C. is a specific type of a finite diagram in $\calC$. It respects the cells and face maps prescribed by the functor $F$ from \eqref{eqn:functorF}. By a reuse of terminology, we say that an object $X \in \calC$ is a G.S.C. if there is a cellular complex $\Phi : J \to \calC$ such that $X = \colim \Phi$. In this case, the functor $\Phi$ will be called a \emph{cellular decomposition} of the object $X$. The integer $L$ is to be represented as the \emph{dimension} of the G.S.C.. The collection of $L$-dimensional geometrical simplicial complexes form a category of its own, which we shall denote as $\SqBrack{ L-\text{complex} }$. The morphisms in this category are induced by graph homomorphisms that preserved the layered structure.

\paragraph{Abstract simplicial complexes} An abstract simplicial complex (A.S.C.) strips away all geometry, coordinates, and ambient spaces. It distills the concept of a "complex" into pure set theory. It is defined as a pair $(V, \Sigma)$, where $V$ is a set of "vertices" and $\Sigma$ is a family of non-empty finite subsets of $V$ (called "abstract simplices"), satisfying closure under subsets. This means that if 
\[\sigma \in \Sigma, \, \tau \subseteq \sigma \imply \tau \in \Sigma.\]
An A.S.C. is just a hypergraph closed under taking subsets. There are no straight lines, no flat triangles, and no topology—just finite sets. Note that A.S.C.-s are independent of $\Topo$ as well as any general category $\calC$ such as ours. However they form an important bridge to convert G.S>C.s into simplicial sets, as we shall see next.

\begin{figure} [!t]
\centering
\begin{tikzpicture}[scale=0.9, transform shape, framed, background rectangle/.style={double, ultra thick, draw=gray, rounded corners}]
	\input{\Diagrams various_complexes.tex}
\end{tikzpicture}
\caption{Various complexes and their relations. See Table \ref{tab:complexes} for more details.}
\label{fig:complexes}
\end{figure}

\paragraph{Some inter-relations} Figure \ref{fig:complexes} and Table \ref{tab:complexes} summarizes these various cell complexes and their inter-relations. Some of these complexes are mutually interchangeable. For example, any geometric simplicial complex $\calK$ carries the information of an abstract simplicial complex. 
\begin{equation} \label{eqn:GSC_to_ASC}
\left\{\begin{array}{c} \mbox{Geometric simplicial} \\ \mbox{complex } \calK \mbox{ with} \\ \mbox{vertex set } V
\end{array} \right\}
\begin{tikzcd}{} \arrow[rr, mapsto] &&{} \end{tikzcd}
\left\{\begin{array}{c} \mbox{Abstract simplicial} \\ \mbox{complex } \calA \mbox{ on } V \,:\, \\ S \in \calA \mbox{ iff } \\ \mbox{vertices span by } S \\ \mbox{is a sub-simplex in } \calK\end{array}\right\}
\end{equation}
In fact, this correspondence is functorial. There is a canonical way in which abstract simplicial complexes can be embedded into the category of simplicial sets. It is created by the correspondence :
\begin{equation} \label{eqn:ASC_to_sSet}
\left\{\begin{array}{c} \mbox{Abstract simplicial} \\ \mbox{complex } \calA \\ \mbox{with ordered} \\ \mbox{vertex set } V \end{array}\right\}
\begin{tikzcd}{} \arrow[rr, mapsto] &&{} \end{tikzcd}
\left\{\begin{array}{c} \mbox{Simplicial set } X : \Delta^{op} \to \SetCat , \\
	X_n := \mbox{all } (n+1)- \mbox{subsets of } V \\ \mbox{which form a face in } \calA \end{array}\right\}
\end{equation}
We have seen in \eqref{eqn:def:Nerve} how the $\Nerve$-functor associates to every $\calC$ object a simplicial set. There is an adjunction 
\[\begin{tikzcd}
\calC \arrow[rrr, bend left=10, "\Nerve"] &&& \sSet \arrow[lll, bend left=10, "\text{Realize}"]
\end{tikzcd}\]
where the geometric realization $\text{Realize}(X)$ of a simplicial set $X: \Delta^{op} \to \SetCat$ is given by the coend
\begin{equation} \label{eqn:def:Realize:coend}
\text{Realize}(X) = \int^{n \in \Delta} X_n \times F_n.
\end{equation}
Here $X_n \times F_n$ is the object formed by taking the disjoint union of copies of $F_n$, one for every element in the set $X_n$. Categorically, this is the copower (or tensor product) of the space $F_n$ by the set $X_n$. A coend is a specific type of colimit that "glues" things together across a category.By definition, the coend $\int^{n \in \Delta} X_n \times F_n$ is constructed by first taking the massive disjoint union (coproduct) of all simplices across all dimensions:
\[\coprod_{n \ge 0} (X_n \times F_n). \]
Then the integral sign quotients this massive space by an equivalence relation generated by the morphisms in $\Delta$. For every morphism $\alpha: [m] \to [n]$ in $\Delta$ (which represents face or degeneracy maps), we enforce the equivalence:
\[(x, \alpha_*(t)) \sim (\alpha^*(x), t), \quad \forall t \in F_m.\]
Here $t \in F_m$ is a point in the topological $m$-simplex, $x \in X_n$ is an abstract $n$-simplex, $\alpha_*: F_m \to F_n$ is the continuous geometric map induced by $\alpha$ (covariant), $\alpha^*: X_n \to X_m$ is the pullback map on the simplicial set induced by $\alpha$ (contravariant).
The realization of a simplicial set as a $\calC$-object can be interpreted equivalently as a right adjoint, as a coend, or as an application of the co-Yoneda lemma. Every presheaf (including any simplicial set $X$) is canonically a colimit of representable functors (the standard simplicial $n$-simplices, $\Delta[n]$). To build any functor $G : \sSet \to \calC$ it is natural to decide where the basic building blocks (the representable functors $\Delta[n]$) should go, and then extend that assignment in a colimit preserving manner to all of $\sSet$.

This completes the presentation of our results and constructions. We end by examining some applications and future directions.

\begin{table}[!t]
\caption{Various complexes used in topology. Most of these complexes can be realized in our abstract categorical framework, as outlined in Section \ref{sec:complexes}.}
\scriptsize
\begin{tabularx}{\textwidth}{|p{0.13\textwidth}|p{0.18\textwidth}|p{0.18\textwidth}|p{0.18\textwidth}|p{0.22\textwidth}|} \hline
	Space Type & Nature of the Space & Chain Group Generators & Computation of Boundary Maps ($\partial_n$) & Interpretation / Primary Application \\ \hline \hline
	Abstract Simplicial Complex &  Combinatorial information in the form of a set system $(V, \Sigma)$ which is closed under subsets. No topology or geometry. & Finite sets of vertices $\{v_0, v_1, \dots, v_n\} \in \Sigma$ spanning an abstract $n$-simplex. & Purely algebraic : $\partial_n[v_0, \dots, v_n]$ equals $\sum_{i=0}^n (-1)^i [v_0, \dots, \hat{v}_i, \dots, v_n]$. It is computed easily via the Smith normal form for matrix reduction. & Provides the computational framework for Topological data analysis (e.g., Vietoris-Rips complexes). The homology of the complex is equivalent to the singular homology of its geometric realization. \\ \hline
	Simplicial Complex & Categorical construct used as an umbrella term for a space equipped with a strict triangulation. & Abstract or geometric $n$-simplices, depending on the context & Combinatorial and determined completely by the underlying abstract vertex scheme. & Key tool for the Simplicial Approximation Theorem. \\ \hline
	Geometric Simplicial Complex & Topological objects embedded in $\mathbb{R}^N$, built from straight rigid simplices intersecting only at proper faces. & Geometric $n$-simplices which are convex hulls of $n+1$ affinely independent points in $\mathbb{R}^N$. & Combinatorial. Boundary maps $\partial_n$ are computed exactly as in the abstract case using the alternating sum of geometric faces. & Used to rigidly compute the Betti numbers and Euler characteristic of triangulable manifolds. \\ \hline
	Cell Complex & A topological space built inductively by gluing topological disks $D^n$ along their boundaries $S^{n-1}$ & Depends on the finiteness. If wild/infinite, falls back to singular $n$-simplices. & Continuous/Singular, depending on the presence of a well-behaved filtration. & A broad classification for spaces built via pushouts. Interprets the broad algebraic invariants of arbitrary glued spaces \\ \hline
	CW Complex & Filtered Topological Space A cell complex satisfying Closure-finiteness and Weak topology & The $n$-cells $e^n_\alpha$. Formally defined as the relative homology $C^{CW}_n(X) = H_n(X^n, X^{n-1}) \cong \mathbb{Z}^{\text{card}(\text{n-cells})}$ & The coefficient $d_{\alpha\beta}$ in $\partial_n(e^n_\alpha) = \sum d_{\alpha\beta} e^{n-1}_\beta$ is the degree of the composite map $S^{n-1} \xrightarrow{\text{attach}} X^{n-1} \xrightarrow{\text{collapse}} S^{n-1}_\beta$ & The most efficient theoretical tool for calculating homology. It radically reduces the number of generators while yielding isomorphic results. \\ \hline
\end{tabularx}
\label{tab:complexes}
\end{table}

\section{Conclusions} \label{sec:conclus}

We have thus seen that under certain assumptions on two distinguished objects $F_0$ and $F_1$ of a general category $\calC$, one gets an inductive procedure by which to build a sequence of cells $F_1, F_1, F_2 , \ldots$. The properties of the pushout present in the wedge construction were leveraged to obtain face and boundary maps between the cells. The main result presented was the formal verification of the fact that these sequence of objects and morphisms assemble to create a co-simplicial object $F : \Delta \to \calC$ satisfying the Axiom of convexity \ref{A:convex}. This functor $F$ along with convexity provides the ingredient for notions of homology, homotopy, and a homotopy invariant notion of homology, as presented in Figure \ref{fig:outline:2}. The entire procedure is summarized in Figure \ref{fig:outline:1}. 
These algebraic invariants provide many of the utilities one sees in the category $\Topo$.  The cells created from $F$ may be simple compared to other objects in $\calC$ and countable in number. However, their strength lie in their ability to be joined in infinitely many ways to create compound objects. There are however some further analogies with $\Topo$ that still remain to be established.

\paragraph{I. Homology calculations} Figure \ref{fig:outline:2} reveals that the homology of any $\calC$-object $X$ is determined solely by the cosimplicial object $F$, via the nerve construction. If there was a way to estimate the homology groups of $X$ externally, such as from measurements, then one has the interesting prospect of being able to reconstruct $X$ from its homology groups. An appealing proposition is to establish \emph{cellular homology} for the various complexes that we have built. Cellular homology is a powerful technique in algebraic topology used to compute the homology groups of a topological space by leveraging its CW complex structure. It is often much more efficient than singular homology because it involves significantly fewer generators. 

Cellular homology builds a chain complex where the "links" are the number of cells in each dimension. The cellular chain group $C_n^{CW}(X)$ is defined as the free abelian group generated by the $n$-cells of $X$, as illustrated below :
\[\dots \xrightarrow{d_{n+1}} C_n^{CW} \xrightarrow{d_n} C_{n-1}^{CW} \xrightarrow{d_{n-1}} \dots \xrightarrow{d_1} C_0^{CW} \xrightarrow{d_0} 0\]
The homology groups are calculated the usual way:
\[H_n^{CW}(X) = \ker(d_n) / \text{im}(d_{n+1}).\]
The utility of cellular homology lies in the boundary map $d_n$. The degree of the attaching map is used to algebraically characterize how an $n$-cell $e^n_\alpha$ maps to an $(n-1)$-cell $e^{n-1}_\beta$. The coefficient $d_{\alpha\beta}$ in the boundary map $d_n(e^n_\alpha) = \sum_\beta d_{\alpha\beta} e^{n-1}_\beta$ is the degree of the map:
\[S^{n-1} \xrightarrow{\phi_\alpha} X^{n-1} \xrightarrow{q_\beta} S^{n-1}, \]
where $\phi_\alpha$ is the attaching map of the $n$-cell and $q_\beta$ collapses everything in the $(n-1)$-skeleton except the interior of the cell $e^{n-1}_\beta$. An interesting proposition to investigate is :

\begin{BigQ}{}{cllr_hmlgy} 
Can cellular homology in this framework be shown to be equivalent to the homology established by the cosimplicial object $F$.
\end{BigQ}

The proof that cellular homology is isomorphic to singular homology for a CW complex is one of the most elegant arguments in algebraic topology. It relies entirely on the long exact sequence of a pair and the excision theorem. This is the next direction in which to expand the results of this article.

\paragraph{II. Excision} In classical topology, if a space $X$ is formed by gluing $A$ and $B$ along their intersection, the following square is a pushout:
\[\begin{array}{ccc}
A \cap B & \longrightarrow & A \\
\downarrow & & \downarrow \\
B & \longrightarrow & X
\end{array}\]
The classical excision theorem implies that applying the singular chain complex functor $C_*$ to this topological pushout square yields a pullback square in the derived category of chain complexes.

\begin{BigQ}{}{excision} 
How can the excision property be realized in our framework ?
\end{BigQ}

There has been several directions in which the excision property was generalized. The most direct abstraction is the concept of an excisive functor \cite{goodwillie1992calculus}, based on the notion of \emph{$\infty$-categories}. A functor $F: \calC \to \calD$ between $\infty$-categories is called \textit{strictly excisive} (or $1$-excisive) if it preserves pushout squares. This means that $F$ maps every homotopy pushout square in $\calC$ to a homotopy pullback square in $\calD$. In this framework, homology is synthetically defined as the \textit{excisive approximation} or the categorical linearization of the unstable identity functor on spaces. 

While singular homology satisfies exact excision, homotopy groups do not. The classical Blakers-Massey theorem or the \emph{Excision theorem for homotopy groups} \cite{massey1991basic} takes a different direction by bounding how far the pushout of a triad is from being a pullback. Recently, this has been generalized \cite{anel2020gen} into a purely categorical theorem valid in any $\infty$-topos. The Generalized Blakers-Massey Theorem categorically quantifies the exact failure of a homotopy pushout square to be a homotopy pullback square, using a chosen modality (factorization system) to measure the connectivity of the maps involved.

In the purely algebraic setting of derived and triangulated categories, excision is abstracted via Verdier localization \cite{verdier1996cat, lurie2017higher}. If $\calC$ is a triangulated category and $\calA$ is a thick subcategory, the Verdier quotient $\calC/\calA$ represents the categorical analogue of relative homology. The "excision isomorphism" in this context states that if $\calA$ and $\calB$ are thick subcategories, then under appropriate conditions:
\[\calB / (\calA \cap \calB) \cong (\calA + \calB) / \calA\]
This perfectly mirrors the topological statement $H_n(B, A \cap B) \cong H_n(X, A)$. Each of these three abstractions of the purely topological phenomenon of excision would be an interesting extension of the present framework.

\paragraph{III. Reconstruction} The idea of finite dimensional complexes in $\calC$ enables a data-driven study of the category $\calC$. Consider the following commutation diagram : 	
\begin{equation} \label{eqn:TDA:scheme}
\begin{tikzcd} [row sep = large]
	\calC \arrow[d, "\text{Homology}"'] && \SqBrack{ L-\text{cellular complex} } \arrow[dashed, Holud]{d}[name=n1]{} \arrow[ll, "\colim"'] \\
	\AbelCat^{\num} \arrow[rr, "\proj_L"'] && \AbelCat^{L} \arrow[bend right=60, Shobuj]{u}[swap, name=n2]{ \text{Reconstruct} }
	\arrow[shorten <=0.1pt, shorten >=4pt, -Bar, to path={(n1) to[out=0,in=180] node[xshift=10pt]{} (n2)} ]{  }
\end{tikzcd}
\end{equation}
The left loop shows the interpretation of $L$ -simplices as $\calC$-objects, the calculation of their Homology, and finally restricting them to the their first $L$ indices. These signature of any $\calC$ object may be directly computable from data. The task of a data-driven reconstruction may be interpreted as the creation of the functor $\text{Reconstruct}$ shown in green, which must be the right adjoint to the composite functor shown in dashed yellow. The interpretation and realization of the reconstruction functor is a subject of future work.

\begin{BigQ}{}{hmlgy_to_smplc} 
Can any finite number of homology groups be realized by a simplicial complex ?
\end{BigQ}

A positive and constructive answer to this question provides a gateway for techniques and numerical methods for topological data analysis \cite[e.g.]{BubenikMilicevic2021, BubenikScott2014prstnt, BauerLesnick2020prsstnc, OtterEtAl2017roadmap} to be applied to the abstract category $\calC$. 

\paragraph{IV. Cellular approximation} The CW Approximation Theorem \cite{Whitehead1949combin, Whitehead1949combin2}  is a cornerstone of algebraic topology, effectively allowing topologists to replace arbitrary spaces with combinatorially manageable ones without changing their weak homotopy type. The theorem establishes that a map between CW complexes which induces isomorphisms on all homotopy groups $\pi_n$ is a true homotopy equivalence.

\begin{BigQ}{}{smplc_to_hmlgy} 
What would be an analog of the cellular approximation theorem in $\calC$ ?
\end{BigQ} 

This question too offers many possible directions of approach, expanding on the usual topological route \cite{Hatcher2002algtopo}; arguments involving  weak homotopy category and adjunctions \cite{may1999algtopo}; or the use of the singular complex and geometric realization functors \cite{Spanier2012algebraic}.

An interesting alternative is Quillen's model categories \cite{Quillen2006hmtpy, hovey1999model} which abstract the formal properties of homotopy theory away from topological spaces. In the standard Strøm or Quillen model structure on $\Topo$, the cofibrant objects are retracts of generalized CW complexes. The CW Approximation Theorem is generalized to the statement that every object in a model category admits a cofibrant replacement. This allowed the exact mechanics of CW approximation to be ported to simplicial sets, chain complexes, and algebraic geometry.

An important result to recall is Milnor's work \cite{milnor1959spaces, LundellWeingram1969topo} showing that the class of spaces with the homotopy type of a CW complex is surprisingly robust. Any space dominated by a CW complex has the homotopy type of a CW complex, and crucially, that the loop space $\Omega X$ of a CW complex has the homotopy type of a CW complex (even though the loop space itself is almost never a CW complex strictly). A generalization of this result would be a stronger statement on the universality of cellular or simplicial complexes

\paragraph{V. General homological categories} Just as homotopical and model categories provide the formal scaffolding to do abstract homotopy theory without topological spaces, a \textit{homological category} \cite{BorceuxBourn2004mal} provides the framework to do homological algebra in strongly non-abelian settings. Classical homological algebra—which gives us tools like exact sequences, the Snake Lemma, and the Nine Lemma—historically required an abelian category (like modules over a ring or sheaves of abelian groups). Abelian categories are highly restrictive because they require morphisms to be added together (enrichment over abelian groups). A homological category  strips away the requirement of commutative addition. It isolates the exact categorical properties needed to manipulate kernels, quotients, and exact sequences in non-commutative algebraic structures. The categorical algebra resulting from homological categories have been well explored \cite[e.g.]{Peschkeetal2024} and several foundational constructions have been discovered \cite[e.g.]{Hartletal2010}.

\paragraph{VI. Triangulation} One of the foundational results in geometry is that every $C^1$-manifold admits a piecewise Linear (PL) triangulation. While initially proved for manifolds of dimensions 2 and 3 \cite{rado1925uber, moise1952affine}, Cairns \cite{cairns1934triangulation}  proved it for submanifolds in Euclidean space, and finally Whitehead  \cite{whitehead1940cmplx} generalized it using the rigorous mechanism of $C^1$-complexes. These papers collectively prove that in dimensions 2 and 3, the distinction between topological and smooth manifolds vanishes; every topological manifold in these dimensions has a unique smooth structure and is therefore triangulable. .
\begin{BigQ}{}{trngltn} 
Characterize the objects of $\calC$ which admit a triangulation, i.e., is isomorphic to a simplicial complex.
\end{BigQ} 

The classification of manifolds has been a long and eventful journey and reflects the theoretical challenges of the abstract challenge posed by this Question. Freedman's classification \cite{freedman1982topo}  of 4-manifolds led to the existence of the $E_8$ manifold, which admits no smooth structure. Shortly after this, Andrew Casson introduced his invariant to definitively prove that such manifolds cannot even admit a combinatorial triangulation). Manolescu \cite{manolescu2016pin} utilized the notion of \textit{Pin(2)-equivariant Seiberg-Witten Floer homology} to show that the Kirby-Siebenmann invariant does not always vanish, proving the existence of non-triangulable topological manifolds in every dimension strictly greater than 4.

We have that seen that if a category $\calC$ is equipped with a terminal object, finite products push-outs has an interval like object $F_1$ one can use that object along with the wedge-construction in an iterative manner to create a co-simplicial object $F$. This object $F$ along with the the interval object $F_1$ create notions of homotopy, homology, and homotopy equivalences which are mutually compatible. Due to this compatibility it becomes meaningful reconstruct or approximate objects in $\calC$ using  simplicial complexes, which could be either a geometric or general cellular complexes.  An important the direction for future work is to show that the resulting cellular homology coincides with the original singular homology. This would lead to a complete cellular approximation theory for the category $\calC$.
\section*{Declarations} 

\paragraph{Acknowledgments} The author is grateful to Dimitri Pavlov (TTU) for his useful feedback.

\paragraph{Funding} This work did not receive any funding or other monetary support.

\paragraph{Conflict of Interest Statement} – The author has no conflicts to disclose.

\paragraph{Ethics Approval Statement} There were no use of any living organisms in this research.

\paragraph{Author Contributions} The sole author SD was responsible for all of the research as well as in the preparation of the document.

\paragraph{Data Availability Statement} The research does not involve any experimental or proprietary data.

\paragraph{Declaration of generative AI use} All text, figures and research in this article was done manually by the author.

\section{Appendix} \label{sec:proofs}
\subsection{Proof of Lemma \ref{lem:Wedge}} \label{sec:proof:Wedge}

To prove that the wedge construct is functorial we have to define its action on morphisms. The diagram below starts with a  morphism $\phi : a\to b$ in $\calC$. 
\begin{equation} \label{eqn:def:Wedge_mrphsm}
\begin{tikzcd} [column sep = large]
	&  & & F_0 \arrow[rrd, Holud] \arrow[rrdd, Akashi] & & \\
	a \arrow[d, "\phi"'] \arrow[r, "="] & a\times F_0 \arrow[rr, "{a\times \tilde{d}_{1,1}}"{pos=0.2}, Holud] \arrow[rru, "!"{pos=0.8}, bend left=49, Holud] \arrow[d, "\phi\times F_0"'] \arrow[rrd, "{\phi \times \tilde{d}_{1,1}}", dashed] & & a\times F_1 \arrow[d, "\phi\times F_1"] \arrow[rr, Holud] & & \Wedge(a) \arrow[Shobuj, d, "\Wedge(\phi)"] \\
	b \arrow[r, "="'] & b\times F_0 \arrow[rr, "{b\times \tilde{d}_{1,1}}"', Akashi] \arrow[rruu, "!"{pos=0.8}, Akashi] & & b\times F_1 \arrow[rr, Akashi] & & \Wedge(b)
\end{tikzcd}
\end{equation}
The yellow and blue sub-diagrams show the separate wedge constructions for the objects $a, b$ respectively. The wedge construct is a pushout and thus universal objects that fulfill a commutation square. Embedded within the above diagram is a co-cone : 
\[\begin{tikzcd} [column sep = large]
& & F_0 \arrow[rrdd, Itranga] & & \\
a\times F_0 \arrow[rr, "{a\times \tilde{d}_{1,1}}"{pos=0.2}, Holud] \arrow[rru, "!"{pos=0.8}, bend left=49, Holud] & & a\times F_1 \arrow[d, "\phi\times F_1", Itranga] & &  \\
& & b\times F_1 \arrow[rr, Itranga] & & \Wedge(b)
\end{tikzcd}\]
The red arrows indicate a co-cone below the yellow diagram. The bottom of this co-cone is $\Wedge(b)$. However $\Wedge(a)$ is the bottom of the universal co-cone below the yellow diagram. Thus the red co-cone must factorize through the universal co-cone. This factorization is via the green arrow shown in \eqref{eqn:def:Wedge_mrphsm}. We can interpret this as the action of $\Wedge$ on the morphism $\phi$. This completes the proof of Lemma \ref{lem:Wedge}. \qed 
\subsection{Identities involving the cone construction} \label{sec:identities}

In this section we explore several identities observed by the morphisms in Diagram \eqref{eqn:conractible:1}.
Since each of the cells are contractible the identity in \eqref{eqn:conractible:2} holds. Substituting the terms there with the redefinitions in \eqref{eqn:Wedge_cells:2} and \eqref{eqn:Wedge_cells:3} we get 
\begin{equation} \label{eqn:conractible:4}
\tilde{s}_{n,n} \circ \tilde{d}_{n+1,n+1} = \Id_{F_n} , \quad \forall n\in\num .
\end{equation}
Next a commutative diagram similar to \eqref{eqn:conractible:1} can be obtained by a second iteration of the wedge construct:
\begin{equation} \label{eqn:double_wedge}
\begin{tikzcd} [column sep = large]
	X \arrow[rrrr, bend left=20, "\Floor_X"] \arrow[rr, Holud, "{X \times \tilde{d}_{1,0}}"] \arrow[d, Itranga, "\Floor_X"'] \arrow[rrdd, Holud] & & X\times F_1 \arrow[d, Itranga, "\Floor_X \times F_1"] \arrow[rr, Holud, "\Pinch_X"] & & \Wedge(X) \arrow[dd, Itranga, "\Wedge(\Floor_X)"] \\
	\Wedge(X) \arrow[drrrr, bend right=30, "\Floor_{\Wedge(X)}"'] \arrow[rr, Akashi, "\Wedge(X) \times \tilde{d}_{1,0}"'] \arrow[rrd, Akashi] & & \Wedge(X) \times F_1 \arrow[rrd, Akashi, "\Pinch_{\Wedge(X)}"] & & \\
	& & 1 \arrow[rruu, Holud] \arrow[rr, Akashi] & & \Wedge^2(X) 
\end{tikzcd}
\end{equation}
The yellow and blue sub-diagrams are the pushout diagrams corresponding to the wedge constructs for $X$ and $\Wedge(X)$ respectively. The vertical arrows in orange represent induced morphisms. The important identity that emerges is
\begin{equation} \label{eqn:wedge_btm_id:1}
\Floor_{\Wedge(X)} \circ \Floor_X = \Wedge \paran{ \Floor_X } \circ \Floor_X, \quad \forall X\in \calC.
\end{equation}
Note that the identity does not require $X$ to be contractible. One similarly obtains the following iterated version of \eqref{eqn:wedge_btm_id:1}
\begin{equation} \label{eqn:wedge_btm_id:2}
\Wedge^{m+1}\paran{ \Floor_{X} } \Floor_{\Wedge^m X} = \Floor_{\Wedge X} \Wedge^{m}\paran{ \Floor_{X} }, \quad \forall X\in \calC,  \,\forall m>0.
\end{equation}
Now suppose that $X$ is contractible. Then the diagram in \eqref{eqn:double_wedge} can be extended as follows :
\[\begin{tikzcd} [column sep = large]
&& && && X \\
X \arrow[rrrr, bend left=20, "\Floor_X"] \arrow[rr, Holud, "{X \times \tilde{d}_{1,0}}"] \arrow[d, Itranga, "\Floor_X"'] \arrow[rrdd, Holud] & & X\times F_1 \arrow[d, Itranga, "\Floor_X \times F_1"] \arrow[rr, Holud, "\Pinch_X"] & & \Wedge(X) \arrow[dd, Itranga, "\Wedge(\Floor_X)"] \arrow[urr, ChhaiC, "\Flatten_{X}"] \\
\Wedge(X) \arrow[drrrr, bend right=30, "\Floor_{\Wedge(X)}"'] \arrow[rr, Akashi, "\Wedge(X) \times \tilde{d}_{1,0}"'] \arrow[rrd, Akashi] & & \Wedge(X) \times F_1 \arrow[rrd, Akashi, "\Pinch_{\Wedge(X)}"] & & \\
& & 1 \arrow[rruu, Holud] \arrow[rr, Akashi] & & \Wedge^2(X) \arrow[drr, ChhaiC, "\Flatten_{\Wedge(X)}"'] \\
&& && && \Wedge(X) 
\end{tikzcd}\]
Next note that the following morphisms may be added while preserving commutations :
\[\begin{tikzcd} [column sep = large]
&& && && X \\
X \arrow[urrrrrr, bend left=20, ChhaiC, "="] \arrow[rrrr, bend left=20, "\Floor_X"] \arrow[rr, Holud, "{X \times \tilde{d}_{1,0}}"] \arrow[d, Itranga, "\Floor_X"'] \arrow[rrdd, Holud] & & X\times F_1 \arrow[d, Itranga, "\Floor_X \times F_1"] \arrow[rr, Holud, "\Pinch_X"] & & \Wedge(X) \arrow[dd, Itranga, "\Wedge(\Floor_X)"] \arrow[urr, ChhaiC, "\Flatten_{X}"] \\
\Wedge(X) \arrow[ddrrrrrr, bend right=30, ChhaiC, pos=0.7, "="'] \arrow[drrrr, bend right=30, "\Floor_{\Wedge(X)}"'] \arrow[rr, Akashi, "\Wedge(X) \times \tilde{d}_{1,0}"'] \arrow[rrd, Akashi] & & \Wedge(X) \times F_1 \arrow[rrd, Akashi, "\Pinch_{\Wedge(X)}"] & & \\
& & 1 \arrow[rruu, Holud] \arrow[rr, Akashi] \arrow[uuurrrr, bend right=10, ChhaiC, pos=0.9, "x"'] \arrow[drrrr, bend right=20, ChhaiC, "\Wedge(x)"'] & & \Wedge^2(X) \arrow[drr, ChhaiC, "\Flatten_{\Wedge(X)}"'] \\
&& && && \Wedge(X) 
\end{tikzcd}\]
Then by uniqueness of pushouts there must be a morphism from $X$ to $\Wedge(X)$ as shown below on as the rightmost dotted arrow :
\[\begin{tikzcd} [column sep = large]
&& && && X \arrow[dotted, dddd, "\Floor_X"] \\
X \arrow[urrrrrr, bend left=20, ChhaiC, "="] \arrow[rrrr, bend left=20, "\Floor_X"] \arrow[rr, Holud, "{X \times \tilde{d}_{1,0}}"] \arrow[d, Itranga, "\Floor_X"'] \arrow[rrdd, Holud] & & X\times F_1 \arrow[d, Itranga, "\Floor_X \times F_1"] \arrow[rr, Holud, "\Pinch_X"] & & \Wedge(X) \arrow[dd, Itranga, "\Wedge(\Floor_X)"] \arrow[urr, ChhaiC, "\Flatten_{X}"] \\
\Wedge(X) \arrow[ddrrrrrr, bend right=30, ChhaiC, pos=0.7, "="'] \arrow[drrrr, bend right=30, "\Floor_{\Wedge(X)}"'] \arrow[rr, Akashi, "\Wedge(X) \times \tilde{d}_{1,0}"'] \arrow[rrd, Akashi] & & \Wedge(X) \times F_1 \arrow[rrd, Akashi, "\Pinch_{\Wedge(X)}"] & & \\
& & 1 \arrow[rruu, Holud] \arrow[rr, Akashi] \arrow[uuurrrr, bend right=10, ChhaiC, pos=0.9, "x"'] \arrow[drrrr, bend right=20, ChhaiC, "\Wedge(x)"'] & & \Wedge^2(X) \arrow[drr, ChhaiC, "\Flatten_{\Wedge(X)}"'] \\
&& && && \Wedge(X) 
\end{tikzcd}\]
The identity \eqref{eqn:conractible:2} then dictates that this morphism is $\Floor_X$. Thus we have derived the commutation 
\begin{equation} \label{eqn:do0p3l:1}
X \mbox{ contractible } \imply 
\begin{tikzcd}
	\Wedge(X) \arrow[d, "\Wedge(\Floor_X)"'] \arrow[rrr, "\Flatten_X"] &&& X \arrow[d, "\Floor_X"] \\
	\Wedge^2(X) \arrow[rrr, "\Flatten_{\Wedge(X)}"'] &&& \Wedge(X)
\end{tikzcd}
\end{equation}
One can draw a similar diagram :
\[\begin{tikzcd} [column sep = large]
&& && && X \arrow[dotted, dddd, "\phi"] \\
X \arrow[d, "\phi"'] \arrow[urrrrrr, bend left=20, ChhaiC, "="] \arrow[rrrr, bend left=20, "\Floor_X"] \arrow[rr, Holud, "{X \times \tilde{d}_{1,0}}"]  \arrow[rrdd, Holud] & & X\times F_1 \arrow[d, Itranga, "\phi \times F_1"] \arrow[rr, Holud, "\Pinch_X"] & & \Wedge(X) \arrow[dd, Itranga, "\Wedge(\phi)"] \arrow[urr, ChhaiC, "\Flatten_{X}"] \\
Y \arrow[ddrrrrrr, bend right=30, ChhaiC, pos=0.7, "="'] \arrow[drrrr, bend right=30, "\Floor_{Y}"'] \arrow[rr, Akashi, "Y \times \tilde{d}_{1,0}"'] \arrow[rrd, Akashi] & & Y \times F_1 \arrow[rrd, Akashi, "\Pinch_{Y}"] & & \\
& & 1 \arrow[rruu, Holud] \arrow[rr, Akashi] \arrow[uuurrrr, bend right=10, ChhaiC, pos=0.9, "x"'] \arrow[drrrr, bend right=20, ChhaiC, "\phi x"'] & & \Wedge(Y) \arrow[drr, ChhaiC, "\Flatten_{\Wedge(X)}"'] \\
&& && && Y
\end{tikzcd}\]
This diagram contains two important commutations
\begin{equation} \label{eqn:do0p3l:2}
X, Y \mbox{ contractible } \imply 
\forall \begin{tikzcd} X \arrow[d, "\phi"'] \\ Y \end{tikzcd} \,:\,
\begin{tikzcd}
	X \arrow[rrrrrr, Akashi, bend left=20, "\cong"] \arrow[rrr, Akashi, "\Floor_X"] \arrow[d, "\phi"'] &&& \Wedge(X) \arrow[d, "\Wedge(\phi)"'] \arrow[rrr, Akashi, "\Flatten_X"] &&& X \arrow[d, "\phi"] \\
	Y \arrow[rrr, Holud, "\Floor_Y"'] \arrow[rrrrrr, Holud, bend right=20, "\cong"'] &&& \Wedge(Y) \arrow[rrr, Holud, "\Flatten_{Y}"'] &&& Y
\end{tikzcd}
\end{equation}

One of the most noticeable features of \eqref{eqn:Wedge_cells:2} is the recursive nature in the definition of the boundary maps $\tilde{d}_{n+1,j}$. 
\[\begin{tikzcd} [column sep = large]
\begin{array}{c} d_{j,j} =\\ \Floor_{F_{j-1}} \end{array} \arrow[r, mapsto, "\Wedge"] \arrow[rr, bend left=20, mapsto, "\Wedge^2"] \arrow[rrrr, bend right=20, mapsto, "\Wedge^m"'] & d_{j+1,j} \arrow[r, mapsto, "\Wedge"] & d_{j+2,j} \arrow[r, mapsto, "\Wedge"] & \cdots \arrow[r, mapsto, "\Wedge"] & d_{j+m,j}
\end{tikzcd}\]
In conclusion, \eqref{eqn:Wedge_cells:2} can be recast as
\begin{equation} \label{eqn:Wedge_cells:4}
d_{j+m,j} = \Wedge^{m} \paran{ \Floor_{F_{j-1}} }, \quad \forall j>0,\, m\geq 0.
\end{equation}
%

\subsection{Proof of simplicial identities} \label{sec:splcl_identities}

\begin{figure}[!t]
\centering
\begin{tikzpicture}[scale=0.9, transform shape]
	\node [style=rect_black_white] (2) at (-0.5\columnA, 0\rowA) {$\begin{array}{l}
			s_{n-1,j} d_{n,i}  \\
			= d_{n-1,i} s_{n-2,j-1} , \\
			0 \leq i \leq n-1 , \\
			0 \leq j\leq n-1 , \\
			i<j
		\end{array}$};
	\node [style=rect_black_white] (2_1) at (1\columnA, 0\rowA) {Induction and Lemma \ref{lem:Wedge}};
	\node [style=rect_black_white] (2_2) at (1\columnA, -1\rowA) {$\begin{array}{l}
			s_{n-1,n-1} d_{n,i}  \\
			= d_{n-1,i} s_{n-2,n-2} , \\
			0 \leq i < n-1
		\end{array}$};
	\node [style=rect_black_white] (2_3) at (1\columnA, 1\rowA) {$\begin{array}{l}
			s_{n-1,n-1} d_{n,n-1}  \\
			= d_{n-1,n-1} s_{n-2,n-2} 
		\end{array}$};
	\node [style=rect_black_white] (2_3_b) at (2\columnA, 1\rowA) {$\begin{array}{l}
			\Flatten\paran{ F_{n-1}, a_{n-1} } \Wedge\paran{ \Floor_{F_{n-2}} }  \\
			= \Floor_{F_{n-2}} \Flatten\paran{ F_{n-2}, a_{n-2} } 
		\end{array}$};
	\node [style=rect_black_white] (2_3_c) at (1.8\columnA, 0\rowA) {Equation \eqref{2_3_b}};
	\node [style=rect_black_white] (2_2_c) at (2.5\columnA, 0\rowA) {Equation \eqref{2_2_b}};
	\node [style=rect_black_white] (2_2_b) at (2\columnA, -1\rowA) {$\begin{array}{l}
			\Flatten\paran{ F_{n-1}, a_{n-1} } d_{n,i}  \\
			= d_{n-1,i} \Flatten\paran{ F_{n-2}, a_{n-2} } , \\
			0 \leq i < n 
		\end{array}$};
	\draw [-to] (2) -- node [above,midway] {$i<n-1, j<n-1$} (2_1);
	\draw [-to] (2) -- node [below,midway] {$i<n-1, j=n-1$} (2_2);
	\draw [-to] (2) -- node [above,midway] {$i=n-1, j=n-1$} (2_3);
	\draw[-to] (2_3) -- node [above,midway] {=} (2_3_b);
	\draw[-to] (2_2) -- node [above,midway] {=} (2_2_b);
	\draw[-to] (2_2_b) -- (2_2_c);
	\draw[-to] (2_3_b) -- (2_3_c);
\end{tikzpicture}
\caption{Proof of Simplicial identity \eqref{eqn:smplc_id}~(ii) in the proof of Theorem \ref{thm:recipe}~(i). The figure shows the case by case analysis of the identity to be proven.}
\label{fig:smplc_id_2}
\end{figure} 

\begin{figure}[!t]
\centering
\begin{tikzpicture}
	\node [style=rect_black_white] (3) at (-0.5\columnA, 0\rowA) {$\begin{array}{l}
			s_{n-1,j} d_{n,j}  \\
			= \Id_{[n-1]} , \\
			0 \leq j\leq n-1 , 
		\end{array}$};
	\node [style=rect_black_white] (3_1) at (1\columnA, 1\rowA) {Induction and Lemma \ref{lem:Wedge}};
	\node [style=rect_black_white] (3_2) at (1\columnA, 0\rowA) {$\begin{array}{l}
			s_{n-1,n-1} d_{n,n-1}  \\
			= \Id_{[n-1]} 
		\end{array}$};
	\node [style=rect_black_white] (3_2_b) at (2\columnA, 0\rowA) {Equation \eqref{eqn:conractible:4}};
	\draw[-to] (3) -- node [above,midway] {$j<n-1$} (3_1);
	\draw[-to] (3) -- node [below,midway] {$j=n-1$} (3_2);
	\draw[-to] (3_2) -- (3_2_b);
\end{tikzpicture}
\caption{Proof of Simplicial identity \eqref{eqn:smplc_id}~(iii) in the proof of Theorem \ref{thm:recipe}~(i). The figure shows the case by case analysis of the identity to be proven.}
\label{fig:smplc_id_3}
\end{figure} 

\begin{figure}[!t]
\centering
\begin{tikzpicture}
	\node [style=rect_black_white] (4) at (-0.5\columnA, 0\rowA) {$\begin{array}{l}
			s_{n-1,j} d_{n,j+1}  \\
			= \Id_{[n-1]} , \\
			0 \leq j\leq n-1 , 
		\end{array}$};
	\node [style=rect_black_white] (4_1) at (1\columnA, 1\rowA) {Induction and Lemma \ref{lem:Wedge}};
	\node [style=rect_black_white] (4_2) at (1\columnA, 0\rowA) {$\begin{array}{l}
			s_{n-1,n-1} d_{n,n}  \\
			= \Id_{[n-1]} 
		\end{array}$};
	\node [style=rect_black_white] (4_2_b) at (2\columnA, 0\rowA) {Equation \eqref{eqn:conractible:4}};
	\draw[-to] (4) -- node [above,midway] {$j<n-1$} (4_1);
	\draw[-to] (4) -- node [below,midway] {$j=n-1$} (4_2);
	\draw[-to] (4_2) -- (4_2_b);
\end{tikzpicture}
\caption{Proof of Simplicial identity \eqref{eqn:smplc_id}~(iv) in the proof of Theorem \ref{thm:recipe}~(i). The figure shows the case by case analysis of the identity to be proven.}
\label{fig:smplc_id_4}
\end{figure} 

\begin{figure}[!t]
\centering
\begin{tikzpicture}[scale=0.9, transform shape]
	\node [style=rect_black_white] (5) at (-0.5\columnA, 0\rowA) {$\begin{array}{l}
			s_{n-1,j} d_{n,i}  \\
			= d_{n-1,i-1} s_{n-2,j} , \\
			0 \leq i \leq n , \\
			0 \leq j\leq n-1 , \\
			i>j+1
		\end{array}$};
	\node [style=rect_black_white] (5_1) at (1\columnA, 1.5\rowA) {Induction and Lemma \ref{lem:Wedge}};
	\node [style=rect_black_white] (5_2) at (1\columnA, 0\rowA) {$\begin{array}{l}
			s_{n-1,j} d_{n,n}  \\
			= d_{n-1,n-1} s_{n-2,j} , \\
			0 \leq j< n-2 
		\end{array}$};
	\node [style=rect_black_white] (5_2_b) at (2\columnA, 0\rowA) {$\begin{array}{l}
			s_{n-1,j} \Floor_{F_{n-1}}  \\
			= \Floor_{F_{n-2}} s_{n-2,j} , \\
			0 \leq j< n-2 
		\end{array}$};
	\node [style=rect_black_white] (5_2_c) at (2\columnA, 1\rowA) {Equation \eqref{eqn:5_2_b}};
	\draw [-to] (2) -- node [above,midway] {$i<n, j<n-1$} (5_1);
	\draw [-to] (5) -- node [above,midway] {$i=n, j<n-1$} (5_2);
	\draw[-to] (5_2) -- node [above,midway] {=} (5_2_b);
	\draw[-to] (5_2_b) -- (5_2_c);
\end{tikzpicture}
\caption{Proof of Simplicial identity \eqref{eqn:smplc_id}~(v) in the proof of Theorem \ref{thm:recipe}~(i). The figure shows the case by case analysis of the identity to be proven.}
\label{fig:smplc_id_5}
\end{figure} 

\begin{figure}[!t]
\centering
\begin{tikzpicture}[scale=0.9, transform shape]
	\node [style=rect_black_white] (6) at (-0.5\columnA, 0\rowA) {$\begin{array}{l}
			s_{n-1,j} s_{n,i}  \\
			= s_{n-1,i} s_{n,j+1} , \\
			0 \leq i \leq n-1 , \\
			0 \leq j\leq n-1 , \\
			i\leq j 
		\end{array}$};
	\node [style=rect_black_white] (6_1) at (1\columnA, 0\rowA) {Induction and Lemma \ref{lem:Wedge}};
	\node [style=rect_black_white] (6_2) at (1\columnA, 1\rowA) {$\begin{array}{l}
			s_{n-1,n-1} s_{n,i}  \\
			= s_{n-1,i} s_{n,n} , \\
			0 \leq i < n-1 
		\end{array}$};
	\node [style=rect_black_white] (6_3) at (1\columnA, -1\rowA) {$\begin{array}{l}
			s_{n-1,n-1} s_{n,n-1}  \\
			= s_{n-1,n-1} s_{n,n} 
		\end{array}$};
	\draw [-to] (6) -- node [above,midway] {$i<n-1, j<n-1$} (6_1);
	\draw [-to] (6) -- node [above,midway] {$i<n-1, j=n-1$} (6_2);
	\draw [-to] (6) -- node [below,midway] {$i=n-1, j=n-1$} (6_3);
\end{tikzpicture}
\caption{Proof of Simplicial identity \eqref{eqn:smplc_id}~(vi) in the proof of Theorem \ref{thm:recipe}~(i). The figure shows the case by case analysis of the identity to be proven.}
\label{fig:smplc_id_6}
\end{figure} 

\begin{figure}[!t]
\centering
\begin{tikzpicture}[scale=0.9, transform shape]
	\node [style=rect_black_white] (1) at (-0.5\columnA, 0\rowA) {$\begin{array}{l}
			d_{n+1,i} d_{n,j}  \\
			= d_{n+1,j+1} d_{n,i} , \\
			0 \leq i \leq n , \\
			0 \leq j\leq n , \\
			i<j
		\end{array}$};
	\node [style=rect_black_white] (1_1) at (1\columnA, 0\rowA) {Induction and Lemma \ref{lem:Wedge}};
	\node [style=rect_black_white] (1_2) at (1\columnA, -1\rowA) {$\begin{array}{l}
			d_{n+1,i} d_{n,n}  \\
			= d_{n+1,n+1} d_{n,i} , \\
			0 \leq i < n 
		\end{array}$};
	\node [style=rect_black_white] (1_3) at (1\columnA, 1\rowA) {$\begin{array}{l}
			d_{n+1,n} d_{n,n}  \\
			= d_{n+1,n+1} d_{n,n} 
		\end{array}$};
	\node [style=rect_black_white] (1_3_b) at (1.9\columnA, 1\rowA) {$\begin{array}{l}
			\Wedge\paran{ \Floor_{F_{n-1}} } \Floor_{F_{n-1}}  \\
			= \Floor_{F_{n}} \Floor_{F_{n-1}} 
		\end{array}$};
	\node [style=rect_black_white] (1_3_c) at (1.8\columnA, 0\rowA) {Equation \eqref{eqn:1_3_c}};
	\node [style=rect_black_white] (1_3_a) at (2.5\columnA, 0\rowA) {Equation \eqref{eqn:1_3_a}};
	\node [style=rect_black_white] (1_2_b) at (2\columnA, -1\rowA) {$\begin{array}{l}
			\Wedge\paran{ d_{n,i} } \Floor_{F_{n-1}} =  \\
			\Floor_{F_{n}} \Wedge\paran{ d_{n-1,i} } , \\
			0 \leq i < n 
		\end{array}$};
	\draw [-to] (1) -- node [above,midway] {$i<n, j<n$} (1_1);
	\draw [-to] (1) -- node [below,midway] {$i<n, j=n$} (1_2);
	\draw [-to] (1) -- node [above,midway] {$i=n, j=n$} (1_3);
	\draw[-to] (1_3) -- node [above,midway] {=} (1_3_b);
	\draw[-to] (1_2) -- node [above,midway] {=} (1_2_b);
	\draw[-to] (1_3_b) -- (1_3_c);
	\draw[-to] (1_2_b) -- (1_3_a);
\end{tikzpicture}
\caption{Proof of Simplicial identity \eqref{eqn:smplc_id}~(i) in the proof of Theorem \ref{thm:recipe}~(i). The figure shows the case by case analysis of the identity to be proven.}
\label{fig:smplc_id_1}
\end{figure}

\paragraph{Proof of \eqref{eqn:1_3_c}} This identity is proved if $X$ in \eqref{eqn:wedge_btm_id:1} is substituted with $F_{n-1}$. 

\paragraph{Proof of \eqref{eqn:1_3_a}} Identities \eqref{eqn:Wedge_cells:4} and \eqref{eqn:wedge_btm_id:2} yield
\[ \Wedge^{n+1-i}\paran{ \Floor_{F_{i-1}} } \Floor_{F_{n-1}} =  d_{n+1,i} d_{n,n}  = d_{n+1,n+1} d_{n,i} = \Floor_{F_{n}} \Wedge^{n-i}\paran{ \Floor_{F_{i-1}} }. \]

\paragraph{Proof of \eqref{2_3_b}} This follows by substituting $X$ in \eqref{eqn:do0p3l:1} by $F_{n-1}$.

\paragraph{Proof of \eqref{2_2_b}} This follows by substituting $X,Y,\phi$ in the right commuting half of \eqref{eqn:do0p3l:2} with $F_{n-2}$, $F_{n-1}$ and $d_{n-1,i}$ respectively.

\paragraph{Proof of \eqref{eqn:5_2_b}} This follows by substituting $X,Y,\phi$ in the left commuting half of \eqref{eqn:do0p3l:2} with $F_{n-2}$, $F_{n-1}$ and $s_{n-2,i}$ respectively.

\subsection{Proof of Lemma \ref{lem:f3sk0}} \label{sec:proof:f3sk0}

To establish the contractibility of $\Wedge(X)$ we need to compare the wedge diagram that creates $\Wedge(X)$, and the wedge diagram by taking a categorical product with $F_1$. For brevity, we shall denote these diagrams as $\angle$ and $\angle\times F_1$ respectively. We begin by observing two natural transformations between $\angle$ and $\angle\times F_1$ : 
\[\begin{tikzcd}[scale cd = 0.6]
& 1 & \\
X \times 1 \times F_1 \arrow[ru, Holud, "!"] \arrow[rrd, Holud, pos=0.7, "X \times \tilde{d}_{1,1} \times F_1"] & & \\
& & X \times F_1 \times F_1 \\
& 1 \arrow[uuu, pos = 0.2, "!"] & \\
X \times 1 \times 1 \arrow[uuu, "{X \times \tilde{d}_{1,0}}"] \arrow[ru, Akashi, "!"] \arrow[rrd, Akashi, "{X \times \tilde{d}_{1,1} \times 1}"'] & & \\
& & X \times F_1 \times 1 \arrow[uuu, "{X \times F_1 \times \tilde{d}_{1,0}}"]
\end{tikzcd} , \quad 
\begin{tikzcd} \Holud{\angle \times F_1} \\ \\ \akashi{\angle} \arrow[uu, Rightarrow, "D_0"] \end{tikzcd} ,\quad
\begin{tikzcd}[scale cd = 0.6]
& 1 & \\
X \times 1 \times F_1 \arrow[ru, Holud, "!"] \arrow[rrd, Holud, pos=0.7, "X \times \tilde{d}_{1,1} \times F_1"] & & \\
& & X \times F_1 \times F_1 \\
& 1 \arrow[uuu, pos = 0.2, "!"] & \\
X \times 1 \times 1 \arrow[uuu, "{X \times  \tilde{d}_{1,1}}"] \arrow[ru, Akashi, "!"] \arrow[rrd, Akashi, "{X \times \tilde{d}_{1,1} \times 1}"'] & & \\
& & X \times F_1 \times 1 \arrow[uuu, "{X \times F_1 \times  \tilde{d}_{1,1}}"]
\end{tikzcd} , \quad 
\begin{tikzcd} \Holud{\angle \times F_1} \\ \\ \akashi{\angle} \arrow[uu, Rightarrow, "D_1"] \end{tikzcd} \]
These natural transformations $D_0, D_1$ embed the diagram $\angle$ within the two ends of $\angle\times F_1$. 
Next note that by assumption we have a homotopy $H$ of $F_1$ into $F_0$, i.e. :
\[\begin{tikzcd} [column sep = large]
& F_1 \times F_1 \arrow[d, Shobuj, "H"] && F_1 \arrow[ll, "F_1 \times \tilde{d}_{1,1}"'] \arrow[dl, "!"] \\
F_1 \arrow[r, "\cong"] \arrow[ur, "F_1 \times \tilde{d}_{1,0}"] & F_1 & 1 \arrow[l]
\end{tikzcd}\]
This $H$ can be used to build a natural transformation $\mathfrak{h}$ between the two wedge-shaped diagrams in the reverse direction : 
\[\begin{tikzcd}[scale cd = 0.7]
& 1 \arrow[ddd, pos=0.2, "="] & \\
X \times 1 \times F_1 \arrow[ddd, "X \times !"'] \arrow[ru, Holud, "!"] \arrow[rrd, Holud, pos=0.7, "X \times \tilde{d}_{1,1} \times F_1"] & & \\
& & X \times F_1 \times F_1 \arrow[ddd, "X \times H"] \\
& 1  & \\
X \times 1 \arrow[ru, Akashi, "!"] \arrow[rrd, Akashi, "{X \times \tilde{d}_{1,1}}"'] & & \\
& & X \times F_1  
\end{tikzcd}\]
If we complete the pushouts of the yellow and blue diagrams respectively we get the following diagram :
\[\begin{tikzcd}[scale cd = 0.7]
& 1 \arrow[drrr, Holud, "\Top_{X \times F_1}"] \arrow[drrrrr, bend left=10, "\Top_{X} \times F_1"] \arrow[ddd, pos=0.2, "="] & \\
X \times 1 \times F_1 \arrow[ddd, "X \times !"'] \arrow[ru, Holud, "!"] \arrow[rrd, Holud, pos=0.7, "X \times \tilde{d}_{1,1} \times F_1"] & & && \Wedge(X\times F_1) \arrow[rr, Shobuj, dotted, "\epsilon"] \arrow[ddd, Shobuj, dotted, pos=0.8, "\mathfrak{h}"] && \Wedge(X) \times F_1 \arrow[dddll, Shobuj, dotted, bend left=20, "G"', "\Wedge(X)\times !"] \\
& & X \times F_1 \times F_1 \arrow[ddd, "X \times H"] \arrow[urr, Holud, "\Floor_{X\times F_1}"'] \arrow[urrrr, bend right=10, pos=0.7, "\Floor_{X}\times F_1"'] \\
& 1\arrow[drrr, bend left=10, Akashi, "\Top_{X \times F_1}"]  & \\
X \times 1 \arrow[ru, Akashi, "!"] \arrow[rrd, Akashi, "{X \times \tilde{d}_{1,1}}"'] & & && \Wedge(X) \\
& & X \times F_1  \arrow[urr, Akashi, bend right=10, "\Floor_{X}"]
\end{tikzcd}\]
The green, dotted arrows in this diagram are created as a consequence of universality of a pushout as a colimit. The morphism $\epsilon$ is induced by the minimality of $\Wedge(X\times F_1)$. The natrual transformation $\mathfrak{h}$ induces a morphism between the two colimits, which has also been denoted as $\mathfrak{h}$ by reuse of notation. Finally, the morphism $G$ which deletes the $F_1$ coordinate commutes with the connecting morphisms of the natural transformation and this commutes with $\epsilon$ and $G$ as shown in the figure. 

If we collect all the diagrams and natural transformation we get the following commutation below on the left.
\[\begin{tikzcd} [scale cd = 0.7]
& \angle \times F_1 \arrow[d, Rightarrow, "\mathfrak{h}"] && \angle \arrow[d, Rightarrow, "!"] \arrow[ll, Rightarrow, "D_1"'] \\
\angle \arrow[ur, Rightarrow, "D_0"] \arrow[r, Rightarrow, "="] & \angle && 1 \arrow[ll, Rightarrow]
\end{tikzcd}
\begin{tikzcd}{} \arrow[rr, mapsto, "\colim"] && {} \end{tikzcd}
\begin{tikzcd} [scale cd = 0.7]
& && & \Wedge(\angle) \arrow[ddl, bend left=10, "!"] \arrow[llld, "D_1"', bend right=10] \\
& \Wedge(\angle \times F_1) \arrow[d, "\mathfrak{h}"] \arrow[rr, "\epsilon"] && \Wedge(\angle) \times F_1  \\
\Wedge(\angle) \arrow[ur, "D_0"] \arrow[r, "="] & \Wedge(\angle) && 1 \arrow[ll]
\end{tikzcd}\]
Applying the colimit functor yields the commutative diagram on the right. Recall that we have established the following commutations so far :
\[\begin{tikzcd} [scale cd = 0.7]
\Wedge(\angle \times F_1) \arrow[d, "\mathfrak{h}"] \arrow[rr, "\epsilon"] && \Wedge(\angle) \times F_1 \arrow[dll, dotted, "\exists G"] \\
\Wedge(\angle)
\end{tikzcd} ,\quad 
\begin{tikzcd} [scale cd = 0.7]
\Wedge(\angle \times F_1) \arrow[rr, "\epsilon"] && \Wedge(\angle) \times F_1  \\
\Wedge(\angle) \arrow[u, "D_0"] \arrow[urr, "\Wedge(\angle) \times \tilde{d}_{1,0}"']
\end{tikzcd},\quad 
\begin{tikzcd} [scale cd = 0.7]
\Wedge(\angle \times F_1) \arrow[rr, "\epsilon"] && \Wedge(\angle) \times F_1  \\
\Wedge(\angle) \arrow[u, "D_1"] \arrow[urr, "\Wedge(\angle) \times \tilde{d}_{1,1}"']
\end{tikzcd}\]
All the morphisms and commutations we have established assemble together to give :
\[\begin{tikzcd}
& \Wedge(\angle) \times F_1 \arrow[d, "G"] && \Wedge(\angle) \arrow[ll, "\Wedge(\angle) \times \tilde{d}_{1,1}"'] \arrow[dl, "!"] \\
\Wedge(\angle) \arrow[r, "\cong"] \arrow[ur, "\Wedge(\angle) \times \tilde{d}_{1,0}"] & \Wedge(\angle) & 1 \arrow[l]
\end{tikzcd}\]
This contraction homotopy, and the fact that $\Wedge(\angle) = \Wedge(X)$ completes the proof of Lemma \ref{lem:f3sk0}. \qed 


\end{document}